\documentclass[12pt,a4paper]{amsart}

\usepackage{amsmath}
\usepackage{amssymb}
\usepackage{amsthm}
\usepackage{latexsym}
\usepackage[T1]{fontenc}
\usepackage{pgf,tikz}
\usetikzlibrary{decorations.pathreplacing}
\usetikzlibrary{arrows}
\usepackage{float}
\usepackage{algorithm}
\usepackage[noend]{algpseudocode}

\theoremstyle{plain}
\newtheorem{theorem}{Theorem}[section]
\newtheorem{corollary}[theorem]{Corollary}
\newtheorem{lemma}[theorem]{Lemma}
\newtheorem{proposition}[theorem]{Proposition}
\newtheorem{conjecture}[theorem]{Conjecture}
\theoremstyle{definition}

\theoremstyle{remark}

\newtheorem{example}{Example}[section]

\numberwithin{equation}{section}

\numberwithin{table}{section}
\def\ds{\displaystyle}
\def\ss{{\mathcal S}}
\def\sgn{\text{sgn}}
\def\GCA{\text{GCA}}
\def\GCB{\text{GCB}}
\def\GCD{\text{GCD}}

\definecolor{ambar}{rgb}{1.0, 0.49, 0.0}
\definecolor{rvwvcq}{rgb}{1,1,1}
\begin{document}

\title[]{Hamilton cycles for involutions of classical types}

\author{Gon\c{c}alo Gutierres}
\address{University of Coimbra, CMUC, Department of Mathematics} \email{ggutc@mat.uc.pt}

\author{Ricardo Mamede}
\address{University of Coimbra, CMUC, Department of Mathematics} \email{mamede@mat.uc.pt}

\author{Jos\'{e} Luis Santos}
\address{University of Coimbra, CMUC, Department of Mathematics} \email{zeluis@mat.uc.pt}

\thanks{Partially supported by the Centre for Mathematics of the University of Coimbra (funded by the Portuguese Government through FCT/MCTES, DOI 10.54499/UIDB/00324/2020)}

%\keywords{Ribbon shapes, Schur functions, Schur support, companion tableau of a Littlewood-Richardson tableau.
%}

%\subjclass[2000]{05A17, 05E05, 05E10, 68Q17}

\maketitle
\begin{abstract}
Let ${\mathcal W}_n$ denote any of the three families of classical Weyl groups: the symmetric groups $\ss_n$, the hyperoctahedral groups (signed permutation groups) $\ss^B_n$, or the
even-signed permutation groups $\ss^D_n$.
In this paper
we give an uniform construction of a Hamilton cycle for the restriction to involutions
on these three families of groups with respect to a inverse-closed connecting set of involutions.
This Hamilton cycle is optimal with respect to the Hamming distance only for the symmetric group $\ss_n$.

We also recall an optimal algorithm for a Gray code for type $B$ involutions.
A modification of this algorithm would provide a Gray Code for  type $D$ involutions with Hamming distance two, which would be optimal. We give such a construction for $\ss^D_4$ and $\ss^D_5$.

%In the previous section we have seen that the Gray Code presented for involutions of type $A$ is optimal with respect to Hamming distance. The same is not true for the codes presented for involutions of types $B$ and $D$. In \cite{grz1} it is proven that the minimal Hamming distance of a Gray Code for type $B$ involutions is two and a Gray Code with Hamming distance two is given. Is this section we will briefly describe that algorithm and, alongside with the description, we will construct the code for $\ss^B_4$ to better illustrate the steps of algorithm.
%A modification of this algorithm would provide a Gray Code for  type $D$ involutions with Hamming distance two, which would be optimal.
%We give such a construction for $\ss^D_4$ and $\ss^D_5$.
%Before starting the description of the code, we introduce some notation for type $B$ involutions. Since $\ss^D_n\subset \ss^B_n$, the some notation is valid for both types.

\end{abstract}


\section{Introduction}

The efficient generation of all objects of a certain class of combinatorial objects is a central problem in enumerative combinatorics, with applications in a vast range of areas ranging from computer science and hardware
or software testing, to biology and biochemistry \cite{mutze, savage}. It can be
used to test hypotheses  or to exhibit combinatorial properties of a class of objects,
 to count the object in a class, and to analyze and prove programs.

A usual approach to this problem has been the generation of the objects of a class as a list in which successive elements differ by a well-defined closeness condition.
An example of such approach is the Binary Reflected Gray Code, or simply Gray code, described and patented in 1953 by Frank Gray \cite{gray}, a researcher at Bell Telephone Laboratories, which generates all $2^n$ $n$-bit strings so that successive strings differ in exactly one bit.
Gray used this code to  prevent spurious output from electromechanical switches, but this code has been widely used in many other areas, such as in circuit design,  data compression and error correction in digital communications.
The theory of Gray codes has
 evolved substantially since F. Gray original work, and the term is now used  in a broader  sense to describe a complete and non-repeating listing of the elements of some class of combinatorial
objects such that successive objects in the listing differ by a well-defined closeness condition.

The concept of a Gray code can be easily translated into graph-theoretical terms. If $\mathcal{W}$ is a class of combinatorial objects, let $G(\mathcal{W})$ be the graph with vertex set $\mathcal{W}$, where two vertices $i$ and $j$ are joined
by an edge whenever $i$ and $j$ satisfy the  closeness condition. The
problem of finding a Gray code for $\mathcal{W}$ is equivalent to find a
Hamiltonian path in $G(\mathcal{W})$. If we demand that the Gray code be closed, that is the initial and final elements in the list must also satisfy the closeness condition, then the problem is equivalent  to the problem of finding a Hamiltonian cycle in $G(\mathcal{W})$.

In this paper we will focus our attention on involutions over the three families of classical Weyl groups: the symmetric groups $\ss_n^A$, the hyperoctahedral groups (signed permutation groups) $\ss^B_n$, and the
even-signed permutation groups $\ss^D_n$. If $\mathcal{W}$ is one of these groups, we consider the Cayley graph $G(\mathcal{W},T)$, for an inverse-closed connecting set $T\subset\mathcal{W}\setminus\{1\}$ (see \cite{meier}). The vertices of this graph are the elements of $\mathcal{W}$, and the edges are given by all possible sets $\{w,w\cdot t\}$, where $w\in\mathcal{W}$ and $t\in T$. The inverse-closed condition of $T$ means that the graph $G(\mathcal{W},T)$  is undirected, since if $\{w,v\}$ is an edge in $G(\mathcal{W},T)$, then $v=w\cdot t$, and also
$w=v\cdot t^{-1}$, so that $\{v,w\}$ is also an edge of $G(\mathcal{W},T)$. Moreover, the graph $G(\mathcal{W},T)$ is connected if and only if $T$ is a generating set for $\mathcal{W}$ \cite{brenti, meier}.

 We will show that the restriction to involutions of the graph
$G(\mathcal{W},T)$ has a Hamiltonian cycle. The choice of the set $T$ is made in order to minimize the maximal distance between two consecutive elements of the cycle.
We adopt as distance between two involutions their Hamming distance, that is, the number of positions in which the words of these two involutions
differ. Using group terminology, the Hamming distance between the involutions $w$ and
$v$ is the number of non-fixed points in the composition $w\cdot v$. The maximal distance
between any two consecutive elements of an Hamilton cycle is the Hamming distance of that
cycle. We will use the Gray codes for types $A$ and $B$ constructed in \cite{grz1} and \cite{grz2} to
generate the Hamilton cycles for these types, and for type $D$ to construct a Gray code  with distance 2 up to $n=5$.
These  Gray codes are proven to be optimal
in relation to the Hamming distance between consecutive elements of the code.

The paper is organized as follows. In Section 2 we describe the classical Weyl groups and its usual combinatorial realizations as the symmetric group, the group of signed permutations, and the even-signed permutation group. Notations and some enumerative results on involutions are addressed for each type. In Section 3 we describe an uniform construction of a Gray code for the involutions in types $A$, $B$ and $D$, based on iterative formulas for their cardinalities,  which translates into  Hamilton cycles in the restriction to involutions of the Cayley graph $G(\ss_n^{\psi},T^{\psi})$,
with $T^{\psi}$ a specific set of reflections of type $\psi\in\{A,B,D\}$. This construction is optimal for type $A$, but not for the other types. Nevertheless, this approach has the advantage of being an uniform construction over the three types, with a simple implementation procedure. In Section 4
we describe a Gray code for type $B$ with optimal Hamming distance, obtained in \cite{grz2},  and generalize this construction to a optimal Gray code for type $D$ involutions. The translation of these codes into Hamilton cycles for the restriction to involutions of the Cayley graph $G(\ss_n^B,T^B)$ and $G(\ss_n^D,T^D)$ are also given.

\section{Involutions in classical Weyl groups}

Let $\mathcal{W}$ be an (irreducible) finite Weyl group with presentation
$$\left\langle s_1,s_2,\ldots,s_n|(s_is_j)^{m_{ij}}=1,\, m_{ii}=1,\text{ and }m_{ij}=m_{ji}\right\rangle.$$
For a Weyl group, the integers $m_{ij}$ take values in the set $\{1,2,3,4,6\}$, and are specified by the Dynkin diagrams represented in Figure \ref{Dynkin}.
The vertices correspond to the generators $s_i$ of $\mathcal{W}$ with  an edge between $i$ and $j$  when $m_{ij}\geq 3$, being that edge  marked with $m_{ij}$ when $m_{ij}\geq 4$.

\begin{center}
\begin{figure}
  \begin{tikzpicture}[scale=.4]
    \draw (-3,0) node[anchor=east]  {$A_n$};
    \foreach \x in {1,...,5}
    \draw[xshift=\x cm,thick] (\x cm,0) circle (.3cm);
    \draw[dotted,thick] (6.25 cm,0) -- +(1.4 cm,0);
    \foreach \y in {1.15,...,2.15,4.15}
    \draw[xshift=\y cm,thick] (\y cm,0) -- +(1.4 cm,0);
    \draw (-1.5,-3) node[anchor=east]  {$B_n\,(n\geq 2)$};
    \foreach \x in {1,...,5}
    \draw[xshift=\x cm,thick] (\x cm,-3) circle (.3cm);
    \draw[dotted,thick] (6.25 cm,-3) -- +(1.4 cm,0);
    \foreach \y in {1.15,...,2.15,4.15}
    \draw[xshift=\y cm,thick] (\y cm,-3) -- +(1.4 cm,0);
    \node at (8.9,-2.4) {4};
    \draw (-1.5,-6) node[anchor=east]  {$D_n\,(n\geq 4)$};
    \foreach \x in {1,...,4}
    \draw[xshift=\x cm,thick] (\x cm,-6) circle (.3cm);
    \draw[xshift=8 cm,yshift=-6cm,thick] (30: 17 mm) circle (.3cm);
    \draw[xshift=8 cm,yshift=-6cm,thick] (-30: 17 mm) circle (.3cm);
    \draw[dotted,thick] (6.25 cm,-6) -- +(1.4 cm,0);
    \foreach \y in {1.15,2.15}
    \draw[xshift=\y cm,thick] (\y cm,-6) -- +(1.4 cm,0);
    \draw[xshift=8 cm,yshift=-6cm,thick] (30: 3 mm) -- (30: 14 mm);
    \draw[xshift=8 cm,yshift=-6cm,thick] (-30: 3 mm) -- (-30: 14 mm);
    \draw (-1.5,-9) node[anchor=east]  {$E_n\, (n=6,7,8)$};
    \foreach \x in {1,...,5}
    \draw[thick,xshift=\x cm,yshift=-10cm] (\x cm,0) circle (3 mm);
    \foreach \y in {1,...,3}
    \draw[thick,xshift=\y cm,yshift=-10cm] (\y cm,0) ++(.3 cm, 0) -- +(14 mm,0);
    \draw[thick] (6 cm,-8 cm) circle (3 mm);
    \draw[thick] (6 cm, -8.3cm) -- +(0, -1.4 cm);
    \draw[dotted,thick] (8.3 cm,-10) -- +(1.4 cm,0);
    \draw (-3,-13) node[anchor=east]  {$F_4$};
    \foreach \x in {1,...,4}
    \draw[xshift=\x cm,yshift=-13cm,thick] (\x cm,0) circle (.3cm);
    \foreach \y in {1,...,3}
    \draw[xshift=\y cm,yshift=-13cm,thick] (\y cm,0) ++(.3 cm, 0) -- +(14 mm,0);
    \node at (5,-12.4) {4};
    \draw (-3,-16) node[anchor=east]  {$G_2$};
    \draw[thick] (2 ,-16) circle (.3 cm);
    \draw[thick] (4 cm,-16) circle (.3 cm);
      \draw[thick] (2.3,-16) -- (3.7,-16);
      \node at (3,-15.4) {6};
  \end{tikzpicture}
  \caption{Dynkin diagrams for irreducible finite Weyl groups.}
  \label{Dynkin}
\end{figure}
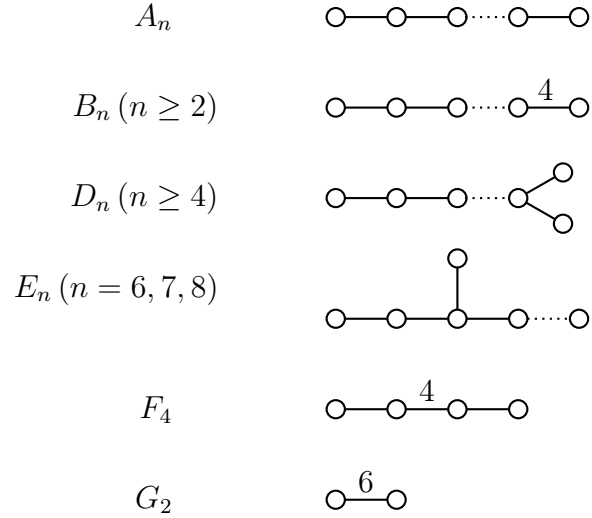
\end{center}

In this paper we focus on the  classical Weyl groups, which correspond to the three infinity families of types $A_{n}$, $B_n$ and $D_n$.
The symmetric group $\ss_{n}$, the hyperoctahedral group $\ss_n^B$, also known as the group of signed permutations, and the even-signed permutation group $\ss_n^D$ are Weyl groups of types $A_{n}$, $B_n$ and $D_n$, respectively.

An involution in $\mathcal{W}$ is an element $w\in\mathcal{W}$ such that $w^2=1$. The elements $ws_iw^{-1}$, with $w\in\mathcal{W}$, are called reflections. Each generator $s_i$ is a reflection, and  any reflection is an involution.
We denote by $I_{n}^A$, $I_n^B$ and $I_n^D$ the set of all involutions of $\ss_{n}^A$, $\ss_n^B$ and $\ss_n^D$, respectively, and by $i^A_n$, $i^B_n$ and $i^D_n$ their cardinality.

\bigskip

The symmetric group $\ss_{n}^A$ is a Weyl group of type $A_{n-1}$, generated by the transpositions of  consecutive integers $s_i^A=(i\;i\!+\!1)$, $i\in[n-1]:=\{1,\ldots,n-1\}$, which satisfy
$$\left(s_i^As_j^A\right)^2=1\,\text{ if }|i-j|>1,\quad\text{ and }\quad \left(s_i^As_{i+1}^A\right)^3=1.$$
The composition of permutations is performed from  left to right.
We use both the decomposition in disjoint cycles and the one-line notation $\pi_1\pi_2\cdots \pi_n$ to represent a permutation  $\pi\in \ss_{n}$, where we set $\pi_i=\pi(i)$ for all $i$. It will be clear from the context
which notation is being used. Multiplying an element $\pi\in\ss_n$ on the right by $s_i$
has the effect of exchanging the values in positions $i$ and $i+1$ in the one-line notation.
The number $i_{n}^A$ of involutions of $\ss_{n}$ satisfies the relation
\begin{equation}\label{eq:recursionA}
i^A_{n}=i^A_{n-1}+(n-1)i^A_{n-2},
\end{equation}
with $i^A_0=i^A_1=1$ (see \cite{Chow}).

\bigskip

Denote by $[\pm n]$ the set $\{\pm 1,\pm 2,\ldots,\pm n\}$ and by $\overline{i}$ the integer $-i$.
The hyperoctahedral group $\ss^B_n$ is the group of all  permutations  $\pi$ of the set
$[\pm n]$ such that
\begin{equation}\label{eq:signperm}
\pi(\,\bar{i}\,)=\overline{\pi(i)},
\end{equation}
for all $i\in[\pm n]$, with composition as group operation, performed from left to right.
Equation \eqref{eq:signperm} indicates that an element of $\ss_n^B$ is entirely defined by its values on $[n]$. For that reason, we will write the elements $\pi\in \ss_n^B$ in one-line notation as $\pi=\pi_1\pi_2\cdots\pi_n$, and call them sign permutations of $[n]$. For example, $\pi=\overline{3}\,2\,\overline{1}\,6\,5\,4\,\overline{7}$ is an element of $\ss_7^B$. Every permutation $\pi\in\ss^B_n$ can be written as a product of disjoint cycles, called the cycle decomposition of $\pi$, obtained by first writing $|\pi|=|\pi_1|\cdots|\pi_n|\in\ss^A_n$ as a disjoint union of cycles, and then placing bars on the letters which have bars in $\pi$. In this setting, if $i,j>0$, the cycle $(i\,j)$ is the permutation that transposes $i$ with $j$, and $\overline{i}$ with $\overline{j}$, the cycle $(\overline{i}\;\overline{j})$ is the permutation that transposes $i$ with $\overline{j}$ and $\overline{i}$ with $j$,  and $(\overline{j})$ is the permutation that transposes $j$ with $\overline{j}$. For instance,
for $\pi=\overline{3}\,2\,\overline{1}\,6\,5\,4\,\overline{7}$ we have $|\pi|=(13)(2)(46)(5)(7)$ and thus the cycle decomposition of $\pi$ is $(\overline{1}\,\overline{3})(2)(4\,6)(5)(\overline{7})$.
Note that the cycle decomposition of a sign involution is a composition of transpositions and cycles of length one. Moreover, using the above notation, the sign of each letter in a transposition $(i\,j)$ with $|i|\neq |j|$ must be the same.

The group $\ss_n^B$ is generated by the involutions
$$s_i^B=(i\,i+1),\text{ for }i\in[n-1],\quad\text{ and }\quad
s_n^B=(\overline{n}),$$
which satisfy
$$(s_i^Bs_{i+1}^B)^3=1\text{ for }i\in[n-2],\quad (s_i^Bs_j^B)^2=1\text{ if }|i-j|>1\quad\text{ and }\quad(s_{n-1}^Bs_n^B)^4=1.$$
Multiplying an element $\pi\in\ss_n^B$ on the right by $s_i^B$, $i\in[n-1]$,
has the effect of exchanging the values in positions $i$ and $i+1$ in the one-line notation. When the generator is $s_n^B$ the multiplication   changes the sign of the value in position $n$. For instance, the permutation $s_5^Bs_6^Bs_7^Bs_6^Bs_5^B$ changes the sign of the letter in position $5$ when applied to an element of $\ss_7^B$:
$$ (\overline{3}\,2\,\overline{1}\,6\,5\,4\,\overline{7})\cdot s_5^Bs_6^Bs_7^Bs_6^Bs_5^B=\overline{3}\,2\,\overline{1}\,6\,\overline{5}\,4\,\overline{7}.$$

The number $i_{n}^B$ of involutions of $\ss_{n}^B$ satisfies the relation
\begin{equation}\label{eq:Bcardformula}
i_{n+1}^B=2i_n^B+2ni_{n-1}^B,
\end{equation}
with $i_0^B=1$ and $i_1^B=2$ (see \cite{Chow}).
\bigskip

The even-signed permutation group $\ss_n^D$ is the subgroup of $\ss_n^B$ consisting of all of the signed permutations
having an even number of negative entries in their one-line notation.
As a set of generators, we take
$$s_i^D=s_i^B=(i\,i+1), \text{ for }i\in[n-1],\quad\text{ and }\quad s_n^D=s_n^Bs_{n-1}^Bs_n^B=(\overline{n-1}\,\,\overline{n}),$$
which satisfy
$$(s_i^Ds_{i+1}^D)^3=(s_{n-2}^Ds_n^D)^3=1\text{ for }i\in[n-2],\quad (s_i^Ds_j^D)^2=(s_{n-1}^Ds_n^D)^2=1\text{ if }|i-j|>1.$$
Multiplying an element $\pi\in\ss_n^D$ on the right by $s_n^D$ has the effect of exchanging the values and the signs in positions $n-1$ and $n$ in the one-line notation. For example,
$$(\overline{3}\,2\,\overline{1}\,6\,\overline{5}\,4\,\overline{7})\cdot s_7^D=\overline{3}\,2\,\overline{1}\,6\,\overline{5}\,7\,\overline{4}.$$

The number of $i_n^D$ of involutions of $\ss_n^D$ can be given by the following recursive formula that uses the number of involutions of type $B$. This formula will allow us to give an uniform construction for the Hamilton cycles in the three types.

\begin{proposition}\label{prop:forrectypeD}
For $n\geq 2$, we have
\begin{equation}\label{eq:Dcardformula}
i_{n+1}^D=i_{n}^B+2n\cdot i_{n-1}^D.
\end{equation}
\end{proposition}
\begin{proof}
Let $\pi=\pi_1\cdots \pi_{n+1}\in\ss_{n+1}^D$ be an involution. Then, either $\pi_{n+1}=\pm(n+1)$ or $\pi_{n+1}=\pm j$, for some $j\in[n]$. In the first case,  $\pi_1\cdots\pi_n$ is an involution in $\ss_n^B$ with $\pi_{n+1}$ positive  if $\pi_1\cdots\pi_n\in \ss_n^D$, and negative otherwise. In the second case, $\pi=(j \,n+1)\cdot \sigma$ or $\pi=(\overline{j}\,\, \overline{n+1})\cdot \sigma$, where $\sigma\in\ss_{n}^D$ with $\sigma_j=j$. The formula now follows from the  bijective correspondence between these involutions $\sigma$ and the involutions of $\ss_{n-1}^D$.
\end{proof}

\section{Recursive Hamilton cycles}

Let $\psi\in\{A,B,D\}$.
In this section we give an uniform construction of an Hamilton cycle for the restriction to involutions of the Cayley graph $G(\ss^{\psi}_n,T^{\psi})$, for a specific generating set $T^{\psi}$ of reflections.
The Hamilton cycles  are achieved by a recursive construction of a Gray code $\text{GC}\psi(n)$ for $I_n^{\psi}$, which satisfy the following two {\bf properties}:
\begin{enumerate}
\item[\bf A1:] The first and last involution of $\text{GC}\psi(n)$ are, respectively, the identity  and the transposition $(n-1\,n)$ or $(\overline{n-1}\,\overline{n})$ (if $\psi\neq A$).
\item[\bf A2:] Any involution in $\text{GC}\psi(n)$ is obtained from its predecessor by a transposition or a rotation of three letters with at most two sign changes (when $\psi\neq A$), or one or two sign changes (when $\psi\neq A$).
\end{enumerate}

To achieve this goal we have to introduce  some notation.

Consider an ordered list $L=(w_1,w_2,\ldots,w_m)$ of involutions over some alphabet $N$ over the integers, which is closed to negation in types $B$ and $D$.
Given letters $i,j,k$ not in $N$, let
$$L\cdot k:=(w_1\cdot k,w_2\cdot k,\ldots, w_m\cdot k)$$
be the ordered list where each involution $w_{\ell}\cdot k$ is the extension of $w_{\ell}$ to the alphabet $N\cup\{k\}$ in type $A$, or $N\cup\{\pm k\}$ otherwise, leaving the letter $k$ fixed. Similarly, define
$$ L\cdot (i\,j):=(w_1\cdot (i\,j),w_2\cdot (i\,j),\ldots, w_m\cdot (i\,j)),$$
where each involution $w_{\ell}\cdot (i\,j)$ is the composition of the involution $w_{\ell}$ with the
transposition $(i\,j)$, now considered over the alphabet $N\cup\{i,j\}$ in type $A$, or  $N\cup\{\pm i, \pm j\}$ otherwise. Running over the elements of $L$ from right to left we get the list
$$\overleftarrow{L}=(w_m,\ldots,w_2,w_1).$$

\begin{example}
If $L=(12, 21)$, then $L\cdot 3=(123,213)$ and $L\cdot (3\,4)=(1243,2143)$.
\end{example}

Let $N=\{a_1,a_2,\ldots,a_n\}$ be an alphabet with $n$ letters, and consider an
involution $\pi\in I_n^A$. Replacing each letter $i$ of $\pi$ by $a_i$ produces an involution $\pi'$ in $N$. Formally, this amounts to defining the bijection $F:[n]\rightarrow N$, where $F(i)=a_i$, which satisfies the composition
$$\pi'=F^{-1}\cdot \pi\cdot F,$$
evaluated from left to right.
Thus, if $\pi=(x_1\,y_1)\cdots(x_k\,y_k)$, then
$\pi'=(a_{x_1}\,a_{y_1})\cdots (a_{x_k}\,a_{y_k})$. To simplify notation, we will write the bijection $F$ in one-line notation $F=a_1a_2\cdots a_n$.

\begin{example}\label{ex:piF}
Consider the involution $\pi=321\in I_3^A$, the alphabet $N=\{1,2,4\}$ and the bijection $F:[3]\rightarrow N$ defined by $F=412$. Then,
$$F^{-1}\cdot \pi \cdot F=142$$
is an involution over the alphabet $N$.
\end{example}

This construction extends naturally to case of sign permutations. For instance, if $\pi=3\overline{2}1\in I_3^B$,  and
$\tilde{F}$ is the extension to type $B$ of the permutation used in Example \ref{ex:piF}, then
$$\tilde{F}^{-1}\cdot \pi \cdot \tilde{F}=\overline{1}42$$
is a sign involution over the alphabet $N$.

Given a list $L=(w_1,w_2,\ldots,w_m)$ of involutions in $I_n^{\psi}$, and alphabet $N$ with $n$ letters and a bijection $F:[n]\rightarrow N$, we denote by
$L^F$ the list
$$L^F=\left(F^{-1}\cdot w_1\cdot F,F^{-1}\cdot w_2\cdot F,\ldots,F^{-1}\cdot w_m\cdot F\right)$$
of involutions over the alphabet $N$. The following result follows from the definitions.

\begin{lemma}\label{lemma:typeAalphabetchange}
Let $N$ be an alphabet with $n$ letters and $F=a_1a_2\cdots a_n$ a bijection between the sets $[n]$ and $N$. If $L$ is a list of involutions in $I_n^{\psi}$, $\psi\in\{A,B,D\}$, satisfying properties {\bf A1} and {\bf A2}, then the sequence $L^F$ is a list of involutions over the alphabet $N$ that also satisfies properties {\bf A1} and {\bf A2}, with $a_{n-1}$ and $a_n$ replaced by $n-1$ and $n$.
\end{lemma}

\subsection{Hamilton cycle of type A}
Consider the Cayley graph  $G(\ss^A_n,T^A)$, where $T^A$ is
formed by reflections and product of two reflections $T^A=X_1^A\cup X_2^A$, where
\begin{align*}
X_1^A&=\{t_{i,j}: i\leq j\}\\
X_2^A&=\{t_{i,j-1}\cdot t_{j,k-1},\, t_{j,k-1}\cdot t_{i,j-1} :i<j<k\},
\end{align*}
  with $t_{i,j}=s_i^As_{i+1}^A\cdots s_{j-1}^As_j^As_{j-1}^A\cdots s_{i+1}^As_i^A$ if $i< j$, and $t_{i,i}=s_i^A$.
Note that a transposition of the letters in positions $i<j$ in a permutation $\pi$ is obtained by the multiplication on the right of $\pi$ by the reflection $t_{i,j-1}$,
while the rotation  of three letters $i<j<k$ of $\pi$ is obtained by multiplication on the right of $\pi$ by $t_{i,j-1}\cdot t_{j,k-1}$ in the case of the rotations $ijk\rightarrow jik\rightarrow jki$, or by $t_{j,k-1}\cdot t_{i,j-1}$ in the case of the rotations
 $ijk\rightarrow ikj\rightarrow kij$.

We can now describe the construction of the Gray code $\GCA(n)$ for the type $A$ involutions of order $n\geq 5$. This construction is achieved by the implementation of the Algorithms \ref{alg:graycodeAodd} and \ref{alg:graycodeAeven} below, one for each parity of $n$, whose recursive calls are triggered by the Gray codes $\GCA(3)$ and $\GCA(4)$, which satisfy properties {\bf A1} and {\bf A2}:
\begin{align}
\GCA(3)&=\big(\text{id},(1\,2),(1\,3),(2\,3)\big)\nonumber\\
\GCA(4)&=\big(\text{id},(1\,3),(1\,3)(2\,4),(2\,4),(1\,4),(1\,4)(2\,3),(2\,3),(1\,2),(1\,2)(3\,4),(3\,4)\big).\label{eq:graycode3}
\end{align}

\begin{algorithm}[H]
\caption{Gray code for the involutions in $\ss_n^A$, $n$ odd ($n\geq 5$)}\label{alg:graycodeAodd}
\begin{algorithmic}[1]

\Procedure{$\GCA(n)$}{}
    \State Set $F=23\cdots(n-1)1$  and write $\GCA^F(n-1)\cdot n$;
    \For{$i=1$ to $\ds\frac{n-1}{2}$}
        \State Set $F=(2i)12\cdots(2i-2)(2i+1)\cdots(n-1)$  and write  $\GCA^F(n-2)\cdot (2i-1\,n)$;
        \State Set $F=(2i-1)12\cdots(2i-2)(2i+1)\cdots(n-1)$  and write $\overleftarrow{\GCA}^F(n-2)\cdot(2i\,n)$;
    \EndFor
\EndProcedure
\end{algorithmic}
\end{algorithm}

%Para i=1, 12...(2i-2) \'{e} vazio

\begin{algorithm}[H]
\caption{Gray code for  involutions in $\ss_n^A$, $n$ even ($n\geq 6$)}\label{alg:graycodeAeven}
\begin{algorithmic}[1]

\Procedure{$\GCA(n)$}{}
    \State Write $\GCA(n-1)\cdot n$;
    \State Set $F=23\cdots(n-1)$  and write $\overleftarrow{\GCA}^F(n-1)\cdot (1\,n)$;
    \For{$i=1$ to $\ds\frac{n}{2}-1$}
        \State Set $F=(2i+1)12\cdots(2i-1)(2i+2)\cdots(n-1)$  and write  $\GCA^F(n-2)\cdot (2i\,n)$;
        \State Set $F=(2i)12\cdots(2i-1)(2i+2)\cdots(n-1)$  and write $\overleftarrow{\GCA}^F(n-2)\cdot(2i+1\,n)$;
    \EndFor
\EndProcedure

\end{algorithmic}
\end{algorithm}

Table \ref{tabela:extipoA} shows the code $\GCA(5)$ generated by Algorithm \ref{alg:graycodeAodd}. The table should be read top to bottom along columns, starting in the leftmost column. Each column contains the involutions generated in the various steps  of the algorithm.

\begin{table}[h]
 \centering
\begin{tabular}{|c|c|c|c|c|}
\hline
12345&52341&15432&12543&21354\\
14325&53241&45312&42513&13254\\
34125&54321&35142&14523&32154\\
32145&52431&15342&21543&12354\\
21345&&&&\\
21435&&&&\\
12435&&&&\\
13245&&&&\\
43215&&&&\\
42315&&&&\\
\hline
\end{tabular}
\caption{The Gray code $\GCA(5)$.}\label{tabela:extipoA}
\end{table}

\begin{theorem}\label{teor:typeAcode}
The application of Algorithm \ref{alg:graycodeAodd} or \ref{alg:graycodeAeven} triggered by Gray codes $\GCA(n-2)$ and $\GCA(n-1)$ satisfying properties {\bf A1} and {\bf A2}, produces a cyclic Gray code for the type $A$ involutions of order $n\geq 5$, which satisfies properties {\bf A1} and {\bf A2}.
\end{theorem}
\begin{proof}
We consider the $n$ odd case, since the even case is analogous.
If $\GCA(n-2)$ and $\GCA(n-1)$ are Gray codes for $I_{n-2}^A$ and $I_{n-1}^A$ satisfying properties {\bf A1} and {\bf A2} then, by Lemma \ref{lemma:typeAalphabetchange}, each involution in any of the sequences  generated in Steps 2, 3 and 4 of Algorithm \ref{alg:graycodeAodd}  is obtained from its predecessor by a transposition or a rotation of three letters.
Note that the leftmost and rightmost involution in each one of these sequences are, respectively:
\begin{itemize}
\item the identity and the transposition $(n\!-\!2\,n\!-\!1)$ in the sequence $\GCA(n-1)^F\cdot n$ obtained in Step 2;
\item  $(2i\!-\!1\,n)$ and $(2i\!-\!1\,n)(n\!-\!2\,n\!-\!1)$ in the sequence $\GCA(n-2)^F\cdot (2i\!-\!1\,n)$ for $i=1,\ldots,(n-1)/2-1$;
\item $(2i\,n)(n\!-\!2\,n\!-\!1)$ and $(2i\,n)$ in the sequence $\overleftarrow{\GCA}(n-2)^F\cdot (2i\!-\!1\,n)$ for $i=1,\ldots,(n-1)/2-1$;
\item $(n\!-\!2\,n)$ and $(n\!-\!2\,n)(n\!-\!4\,n\!-\!3)$ in the sequence  $\GCA(n-2)^F\cdot (2i\!-\!1\,n)$ for $i=(n-1)/2$;
\item $(n\!-\!1\,n)(n\!-\!4\,n\!-\!3)$ and $(n\!-\!1\,n)$ in the sequence $\overleftarrow{\GCA}(n-2)^F\cdot (2i\!-\!1\,n)$ for $i=(n-1)/2$.
\end{itemize}
It follows that $\GCA(n)$ is a cyclic Gray code for $I_n^A$ satisfying properties
{\bf A1} and {\bf A2}.

Finally, notice that by its construction, all involutions generated by the algorithm are pairwise distinct.  Moreover, the number of involutions in all these sequences is equal to $t_{n-1}^A+(n-1)t_{n-2}^A$, which gives the total amount of involutions in $I_n^A$ by equation \eqref{eq:recursionA}.
\end{proof}

The Gray code $\GCA(n)$ translates into an Hamilton cycle in the restriction to involutions of the Cayley graph $G(\ss_n^A,T^A)$.

\begin{corollary}
The Gray code $\GCA(n)$ produced by Algorithm \ref{alg:graycodeAodd} or \ref{alg:graycodeAeven} is an Hamilton cycle in the restriction to involutions of the Cayley graph $G(\ss_n^A,T^A)$, with Hamming distance three, for $n\geq 3$. This is the minimal Hamming distance of any Gray code for $I_n^A$ with $n\geq 3$.
\end{corollary}
\begin{proof}
Since $\GCA(n)$ satisfies properties {\bf A1} and {\bf A2}, each involution in the list is obtained from its predecessor by a transposition or a rotation of three letters.
  Since a transposition of two letters in an involution $\pi$ is obtained by the multiplication on the right of $\pi$ by the reflection $t_{i,j}$, and the rotation of three letters is obtained by multiplication on the right of $\pi$ by the product of two reflections $t_{i,j-1}\cdot t_{j,k-1}$ or $t_{j,k-1}\cdot t_{i,j-1}$, it follows that any two consecutive involutions in the code $\GCA(n)$ are connected by an edge in the Cayley graph $G(\ss_n^A,T^A)$. The same is true for the initial and final elements of the code. Thus, $\GCA(n)$ is an Hamilton cycle in the restriction to involutions of the Cayley graph $G(\ss_n^A,T^A)$.

  It follows that the Hamming distance of this code is three. We will prove that three is the minimal distance in a code for $I_n^A$. The distance between two involutions is at least two, and this number is achieved when the involutions differ by exactly one transposition. Thus, a code for $I_n^A$ with distance two must be a sequence of involutions that differ by exactly one transposition, that is an alternating  sequence of even and odd involutions. But since the excess of even to odd involutions is greater than $1$ for $n>2$ (see \cite{excessevenodd}), there cannot be such sequence. Therefore, we conclude that three is the minimal distance for a code for $I_n^A$.
\end{proof}

\subsection{Hamilton cycle of type B}
The Cayley graph $G(\ss^B_n,T^B)$ for the involutions of type $B$ are generated by the set of reflections and product of two reflections $T^B=X_1^B\cup X_2^B\cup X_3^B\cup X_4^B$, where
\begin{align*}
X_1^B&=\{t_{i,j}: i\leq j\}\\
X_2^B&=\{t_{i,n}\cdot t_{j,n}:i<j<n\}\\
X_3^B&=\{t_{i,j}\cdot t_{j+1,n} :i<j<n\}\\
X_4^B&=\{t_{j,k-1}\cdot t_{i,j-1}:i<j<k\},
\end{align*}
with $t_{i,j}=s_i^Bs_{i+1}^B\cdots s_{j-1}^Bs_j^Bs_{j-1}^B\cdots s_{i+1}^Bs_i^B$ if $i\leq j$, and $t_{i,i}=s_i^B$.
The construction of the Gray code $\GCB(n)$ for the involutions  of order $n\geq 4$  in $\ss^B_n$ is obtained implementing Algorithm \ref{alg:typeB}, whose recursive calls are triggered by the Gray codes
$$\GCB(2)=\big(id, (\overline{1}), (\overline{1})(\overline{2}), (\overline{2}), (\overline{1}\,\overline{2}), (12) \big)$$
and
\begin{align*}
\GCB(3)=&\big(id, (\overline{1}), (\overline{1})(\overline{2}), (\overline{2}), (\overline{1}\,\overline{2}), (12),
(12)(\overline{3}), (\overline{1}\,\overline{2})(\overline{3}), (\overline{2})(\overline{3}), (\overline{1})(\overline{2})(\overline{3}),(\overline{1})(\overline{3}), (\overline{3}),\\
&(\overline{1}\,\overline{3}), (\overline{1}\,\overline{3})(\overline{2}), (13)(\overline{2}), (13), (23), (\overline{1})(23), (\overline{1})(\overline{2}\,\overline{3}),
(\overline{2}\,\overline{3})\big)
\end{align*}
which satisfy properties {\bf A1} and {\bf A2}.

\begin{algorithm}[h]
\caption{Gray code for the involutions in $\ss_n^B$ ($n\geq 4$)}\label{alg:typeB}
\begin{algorithmic}[1]

\Procedure{$\GCB(n)$}{}
    \State Write $\GCB(n-1)\cdot n$;
    \State Write $\overleftarrow{\GCB}(n-1)\cdot \overline{n}$;
    \For{$i=1$ to $n-1$}
        \State Set $F=12\cdots(i-1)(i+1)\cdots(n-1)$;
        \If{$i$ is odd}
            \State Write  $\GCB^F(n-2)\cdot (\overline{i}\,\overline{n}), \overleftarrow{\GCB}^F(n-2)\cdot(i\,n)$;
        \Else
            \State Write $\GCB^F(n-2)\cdot(i\,n),\overleftarrow{\GCB}^F(n-2)\cdot (\overline{i}\,\overline{n})$;
        \EndIf
    \EndFor
\EndProcedure

\end{algorithmic}
\end{algorithm}

Table   \ref{tabela:extipoB} shows the Gray code $\GCB(4)$ obtained by Algorithm \ref{alg:typeB}, which should be read down columns, from left to right.
Each column contains the involutions generated in a different step of the algorithm.

\begin{table}[h]
 \centering
\begin{tabular}{|c|c|c|c|c|c|}
\hline
  $1\,2\,3\,4$&                                    $1\,\overline{3}\,\overline{2}\,\overline{4}$&                    $\overline{4}\,2\,3\,\overline{1}$& $1\,4\,3\,2$ & $1\,2\,\overline{4}\,\overline{3}$     \\
 $\overline{1}\,2\,3\,4$&                           $\overline{1}\,\overline{3}\,\overline{2}\,\overline{4}$&        $\overline{4}\,\overline{2}\, 3\,\overline{1}$ &$\overline{1}\,4\,3\,2$ &$\overline{1}\,2\,\overline{4}\,\overline{3}$     \\
 $\overline{1}\,\overline{2}\,3\,4$&                $\overline{1}\, 3\, 2\,\overline{4}$&                            $\overline{4}\,\overline{2}\,\overline{3}\,\overline{1}$ &$\overline{1}\,4\,\overline{3}\,2$ & $\overline{1}\,\overline{2}\,\overline{4}\,\overline{3}$     \\
  $1\,\overline{2}\,3\,4$&                             $1\, 3\, 2\,\overline{4}$&                                      $\overline{4}\,2\,\overline{3}\,\overline{1}$ & $1\,4\,\overline{3}\,2$ & $1\,\overline{2}\,\overline{4}\,\overline{3}$  \\
 $\overline{2}\,\overline{1}\,3\,4$&                 $3\, 2\, 1\,\overline{4}$&                                      $\overline{4}\,\overline{3}\,\overline{2}\,\overline{1}$ &$\overline{3}\,4\,\overline{1}\,2$ & $\overline{2}\,\overline{1}\,\overline{4}\,\overline{3}$     \\
  $2\,1\,3\,4$&                                      $3\,\overline{2}\, 1\,\overline{4}$&                            $\overline{4}\,3\,2\,\overline{1}$& $3\,4\,1\,2$ & $2\,1\,\overline{4}\,\overline{3}$ \\
  $2\,1\,\overline{3}\,4$&                          $\overline{3}\,\overline{2}\,\overline{1}\,\overline{4}$&         $4\,3\,2\,1$& $3\,\overline{4}\,1\,\overline{2}$&$2\,1\,4\,3$\\
 $\overline{2}\,\overline{1}\,\overline{3}\,4$&     $\overline{3}\, 2\,\overline{1}\,\overline{4}$&                   $4\overline{3}\,\overline{2}\,1$&$\overline{3}\,\overline{4}\,\overline{1}\,\overline{2}$&$\overline{2}\,\overline{1}\,4\,$3\\
  $1\,\overline{2}\,\overline{3}\,4$&                $1\, 2\,\overline{3}\,\overline{4}$&                             $4\,2\,\overline{3}\,1$& $1\,\overline{4}\,\overline{3}\,\overline{2}$&$1\,\overline{2}\,4\,3$\\
 $\overline{1}\,\overline{2}\,\overline{3}\,4$&     $\overline{1}\, 2\,\overline{3}\,\overline{4}$&                   $4\,\overline{2}\,\overline{3}\,1$&$\overline{1}\,\overline{4}\,\overline{3}\,\overline{2}$& $\overline{1}\,\overline{2}\,4\,3$\\
 $\overline{1}\,2\,\overline{3}\,4$&                $\overline{1}\,\overline{2}\,\overline{3}\,\overline{4}$&         $4\,\overline{2}\,3\,1$&$\overline{1}\,\overline{4}\,3\,\overline{2}$&$\overline{1}\,2\,4\,3$\\
  $1\,2\,\overline{3}\,4$&                             $1\,\overline{2}\,\overline{3}\,\overline{4}$&                   $4\,2\,3\,1$& $1\,\overline{4}\,3\,\overline{2}$& $1\,2\,4\,3$\\
 $\overline{3}\,2\,\overline{1}\,4$&                $\overline{2}\,\overline{1}\,\overline{3}\,\overline{4}$&    &&\\
 $\overline{3}\,\overline{2}\,\overline{1}\,4$&      $2\, 1\,\overline{3}\,\overline{4}$&                        &&\\
  $3\,\overline{2}\,1\,4$&                           $2\, 1\, 3\,\overline{4}$&                                  &&\\
  $3\,2\,1\,4$&                                     $\overline{2}\,\overline{1}\, 3\,\overline{4}$&              &&\\
  $1\,3\,2\,4$&                                      $1\,\overline{2}\, 3\,\overline{4}$&                        &&\\
 $\overline{1}\,3\,2\,4$&                           $\overline{1}\,\overline{2}\, 3\,\overline{4}$&              &&\\
 $\overline{1}\,\overline{3}\,\overline{2}\,4$&     $\overline{1}\, 2\, 3\,\overline{4}$&                        &&\\
  $1\,\overline{3}\,\overline{2}\,4$&                $1\, 2\, 3\,\overline{4}$&                                  &&\\
\hline
\end{tabular}
\caption{The Gray code $\GCB(4)$.}\label{tabela:extipoB}
\end{table}

\begin{theorem}\label{teor:typeBcode}
The application of Algorithm \ref{alg:typeB} with input the Gray codes $\GCB(n-2)$ and $\GCB(n-1)$, which satisfy properties {\bf A1} and {\bf A2}, produces a cyclic Gray code for the type $B$ involutions of order $n\geq 4$, which satisfies properties {\bf A1} and {\bf A2}.
\end{theorem}
\begin{proof}
As in the proof of Theorem \ref{teor:typeAcode} for the type $A$ case, now using the recursive formula \eqref{eq:Bcardformula} for the number of type $B$ involutions, we can show that $\GCB(n)$ is a list of all involutions in $\ss_n^B$. Assume that $\GCB(n-2)$ and $\GCB(n-1)$
satisfies properties {\bf A1} and {\bf A2}. Moreover, assume that whenever two consecutive elements of one of these codes differ by a rotation of three letters, there is no sign change, and whenever there is a transposition of two letters, at most one of these letters changes its sign. Notice that this is certainly true for  the trigger cases $\GCB(2)$ and $\GCB(3)$. We will show that the same is true for the sequence $\GCB(n)$ generated by  Algorithm \ref{alg:typeB}. A rotation of three letters between two consecutive involutions in $\GCB(n)$ only occurs between two  iterations of the cycle {\bf for} in Step 4 of the algorithm, with no sign changes. Similarly,  when a transposition occur between two consecutive involutions, there is at most one sign change. This situation only occurs   between the  involutions $(\overline{n})$ and  $(\overline{1}\,\overline{n})$, generated in Step 3 and  in the first iteration of the cycle {\bf for} in Step 4, respectively, with the letter 1 changing  its sign. By Lemma \ref{lemma:typeAalphabetchange} we conclude that $\GCB(n)$ satisfy properties {\bf A1} and {\bf A2}, and whenever two consecutive letters of one of these codes differ by a rotation of three letters, there is no sign change, and whenever there is a transposition of two letters, at most one of these letters changes its sign.
\end{proof}

As in the type $A$ case,
the Gray code $\GCB(n)$ translates into an Hamilton cycle in the restriction to involutions of the Cayley graph $G(\ss_n^B,T^B)$ with Hamming distance three, for $n\geq 3$. However, in this case, this distance is not the minimal possible distance for such a Gray code, as it will be discussed in Section \ref{sec:optimal_algB}.

\begin{corollary}\label{cor:hamtypeB}
The Gray code $\GCB(n)$ produced by Algorithm \ref{alg:typeB} is an Hamilton cycle in the restriction to involutions of the Cayley graph $G(\ss_n^B,T^B)$, with Hamming distance three, for $n\geq 4$.
\end{corollary}
\begin{proof}
Since $\GCB(n)$ satisfies properties {\bf A1} and {\bf A2}, each involution in the list is obtained from its predecessor by a transposition or a rotation of three letters and at most two sign changes. Moreover, by the proof of Theorem \ref{teor:typeBcode} a rotation of three letters only occur  with no sign changes. Additionally,  when a transposition between letters occurs, there is at most one sign change on the smallest of the two letters.  Since a sign change in position $i$ of a involution $\pi$ is obtained by the multiplication on the right of $\pi$ by $t_{i,n}$, %a transposition of two letters in positions $i<j+1$ in $\pi$ is obtained by the multiplication on the right  by the reflection $t_{i,j}$,
a transposition of two letters in positions $i<j+1$ in $\pi$ followed by a sign change of the letter in position $j+1$ is obtained by the multiplication on the right of $\pi$  by  $t_{i,j}\cdot t_{j+1,n}$, and the rotation of three letters in positions $i<i+1<j$ is obtained by multiplication on the right of $\pi$ by the product of two reflections $t_{i+1,j-1}\cdot t_{i,i}$, it follows that any two consecutive involutions in the code $\GCB(n)$ have Hamming distance at most 3 and are connected by an edge in the Cayley graph $G(\ss_n^B,T^B)$. The same is true for the initial and final elements of the code. Thus, $\GCB(n)$ is an Hamilton cycle in the restriction to involutions of the Cayley graph $G(\ss_n^B,T^B)$, and  their Hamming distance is three.
  \end{proof}

\subsection{Hamilton cycle of type D}\label{section:typeD} To generate an Hamilton cycle for the involutions in $\ss_n^D$, consider the generating set
$T^D=X_1^D\cup X_2^D\cup X_3^D\cup X_4^D$, where
\begin{align*}
X_1^D&=\{t_{i,j}: i\leq j<n\}\\
X_2^D&=\{t_{i,n}\cdot t_{j,n}:i<j<n\}\\
X_3^D&=\{t_{i,j}\cdot t_{j+1,n}\cdot t_{k,n}:i<j<k-1\}\\
X_4^D&=\{t_{j,k-1}\cdot t_{i,j-1},\,\, t_{j,k-1}\cdot t_{i,j-1}\cdot t_{i,n}\cdot t_{j,n},\,\, t_{j,k-1}\cdot t_{i,j-1}\cdot t_{i,n}\cdot t_{k,n} :i<j<k\},
\end{align*}
with $t_{i,j}=s_i^Ds_{i+1}^D\cdots s_{j-1}^Ds_j^Ds_{j-1}^D\cdots s_{i+1}^Ds_i^D$ if $i< j<n$, and  $t_{i,i}=s_i^D$ for $i<n$,
$$t_{n-1,n}=s_{n-1}^Ds_n^D\quad\text{ and }\quad t_{i,n}=s_i^Ds_{i+1}^D\cdots s_{n-2}^D\cdot\left(s_{n-1}^Ds_n^D\right)\cdot s_{n-2}^D\cdots s_{i+1}^Ds_i^D,$$
for $i\in[n-2]$.

Given a list $L=(w_1,w_2,\ldots,w_m)$ of involutions of $\ss_{n-1}^B$ we write $L\cdot\tilde{n}=(w_1\cdot\tilde{n},w_2\cdot\tilde{n},\ldots,w_m\cdot\tilde{n})$,
where each  $w_i\cdot\tilde{n}$ is $w_i\cdot n$ if the involution $w_i$ has an even number of negative signs, and is
$w_i\cdot \overline{n}$ otherwise. Each permutation in $L\cdot\tilde{n}$ is an involution in $\ss_n^D$.

The construction of the Gray code $\GCD(n)$ for the involutions  of order $n\geq 4$  in $G(\ss^D_n,T^D)$ is obtained implementing Algorithm \ref{alg:typeD}, whose recursive calls are triggered by the type $B$ Gray codes $\GCB(n)$ and the codes
$$\GCD(1)=(id)\quad\text{ and }\quad\GCD(2)=\left(id, (\overline{1})(\overline{2}), (\overline{1}\,\overline{2}), (1\,2)\right).$$

\begin{algorithm}[H]
\caption{Gray code for the involutions in $\ss_n^D$ ($n\geq 3$)}\label{alg:typeD}
\begin{algorithmic}[1]

\Procedure{$\GCD(n)$}{}
    \State Set $F=23\cdots(n-1)1$ and write $\GCB^F(n-1)\cdot \widetilde{n}$;
    \For{$i=1$ to $n-1$}
        \State Set $F=12\cdots(i-1)(i+1)\cdots(n-1)$;
        \State Write  $\GCD^F(n-2)\cdot (i\,n)$;
        \State Write $\overleftarrow{\GCD}^F(n-2)\cdot(\overline{i}\,\overline{n})$;
    \EndFor
\EndProcedure

\end{algorithmic}
\end{algorithm}

Tables   \ref{tabela:extipoD3} and \ref{tabela:extipoD4}   show the Gray codes $\GCD(3)$ and $\GCD(4)$ obtained by Algorithm \ref{alg:typeB}, which as before,  should be read down columns, from left to right.
Each column contains the involutions generated in a different step of the algorithm.

\begin{table}[h]
 \centering
\begin{tabular}{|c|c|c|}
\hline
$1\,2\,3$                      &    $3\,2\,1$&$1\,3\,2$\\
$1\,\overline{2}\,\overline{3}$&$\overline{3}\,2\,\overline{1}$&$1\,\overline{3}\,\overline{2}$\\
$\overline{1}\,\overline{2}\,3$&   &     \\
$\overline{1}\,2\,\overline{3}$&    &    \\
$\overline{2}\,\overline{1}\,3$&    &    \\
$2\,1\,3$&&\\
\hline
\end{tabular}
\caption{The Gray code $\GCD(3)$.}\label{tabela:extipoD3}
\end{table}

\begin{table}[h]
 \centering
\begin{tabular}{|c|c|c|c|}
\hline
     $ 1\,          2\,      3\,     4$&                              $4\,     2\,     3\,     1$                               &  $1\,     4\,     3\,     2$&                                 $ 1\,     2\,     4\,     3$  \\
     $ 1\,\overline{2}\,          3\,    \overline{4}$&               $4\,\overline{2}\,\overline{3}\,1$                        &  $\overline{1}\,4\,\overline{3}\,2$&                          $\overline{1}\,\overline{2}\,4\,     3$  \\
     $ 1\,\overline{2}\,\overline{3}\,     4$&                        $4\,\overline{3}\,\overline{2}\,1$                        &  $\overline{3}\,4\,\overline{1}\,2$&                          $\overline{2}\,\overline{1}\,4\,3$  \\
     $ 1\,     2\,\overline{3}\,\overline{4}$&                        $4\,      3\,     2\,     1$                              &  $3\,     4\,     1\,     2$&                                 $ 2\,     1\,     4\,     3$  \\
     $ 1\,\overline{3}\,\overline{2}\,     4$&                       $\overline{4}\,3\,2\,\overline{1}$                         &  $3\,\overline{4}\,1\,\overline{2}$&                          $ 2\,1\,\overline{4}\,\overline{3}$  \\
     $ 1\,     3\,     2\,     4$&                                   $\overline{4}\,\overline{3}\,\overline{2}\,\overline{1}$   &  $\overline{3}\,\overline{4}\,\overline{1}\,\overline{2}$&    $\overline{2}\,\overline{1}\,\overline{4}\,\overline{3}$   \\
     $\overline{1}\,     3\,     2\,    \overline{4}$&               $\overline{4}\,\overline{2}\,\overline{3}\,\overline{1}$   & $\overline{1}\,\overline{4}\,\overline{3}\,\overline{2}$&     $\overline{1}\,\overline{2}\,\overline{4}\,\overline{3}$  \\
     $\overline{1}\,\overline{3}\,\overline{2}\,\overline{4}$&       $\overline{4}\,2\,3\,\overline{1}$                         &  $1\,\overline{4}\,3\,\overline{2}$&                          $ 1\,2\,\overline{4}\,\overline{3}$    \\
     $\overline{1}\,     2\,\overline{3}\,     4$&                  & &   \\
     $\overline{1}\,\overline{2}\,\overline{3}\,\overline{4}$&      & &    \\
     $\overline{1}\,\overline{2}\,     3\,     4$&                  & &     \\
     $\overline{1}\,     2\,     3\,\overline{4}$&                  & &      \\
     $\overline{2}\,\overline{1}\,     3\,     4$&                  & &       \\
     $\overline{2}\,\overline{1}\,\overline{3}\,\overline{4}$&      & &        \\
     $ 2\,     1\,\overline{3}\,\overline{4}$&                      & &     \\
     $ 2\,     1\,     3\,     4$&                                  & &      \\
     $ 3\,     2\,     1\,     4$&                                  & &       \\
     $ 3\,\overline{2}\,     1\,\overline{4}$&                      & &        \\
     $\overline{3}\,\overline{2}\,\overline{1}\,\overline{4}$&      & &         \\
     $\overline{3}\,     2\,\overline{1}\,     4$&                  & &          \\

\hline
\end{tabular}
\caption{The Gray code $\GCD(4)$.}\label{tabela:extipoD4}
\end{table}

\begin{theorem}\label{teor:algtypeD}
If the Gray codes $\GCD(n-2)$ and $\GCB(n-1)$ satisfy properties {\bf A1} and {\bf A2}, then
the Algorithm \ref{alg:typeD}  produces a cyclic Gray code $\GCD(n)$ for the type $D$ involutions of order $n\geq 3$, which satisfies properties {\bf A1} and {\bf A2}.
\end{theorem}
\begin{proof}
As in the proof of  Theorem \ref{teor:typeAcode} for the type $A$ case, now using the formula \eqref{eq:Dcardformula} for the number of type $D$ involutions, obtained in Proposition \ref{prop:forrectypeD}, we can show that the sequence $\GCD(n)$ given by Algorithm \ref{alg:typeD} contains all involutions in $\ss_n^D$, each appearing exactly once.

If $\GCD(n-2)$ and $\GCB(n-1)$ satisfy properties {\bf A1} and {\bf A2}, then by construction  the sequences  $\GCB^F(n-1)\cdot \widetilde{n}$, $\GCD^F(n-2)\cdot (i\,n)$ and $\overleftarrow{\GCD}^F(n-2)\cdot(\overline{i}\,\overline{n})$ generated in Steps 2, 5 and 6 of the algorithm, are formed by pairwise distinct involutions and also satisfy property {\bf A2}. It remains to check that the involutions connecting these sequences also satisfy property {\bf A2}. The last element of $\GCB^F(n-1)\cdot \widetilde{n}$ is the involution $$u_1=(n-1)\,2\cdots\,1\,n\quad\text{ or }\quad u_2=(\overline{n-1})\,2\cdots\,\overline{1}\,n,$$
and the first element of $\GCD^F(n-2)\cdot (i\,n)$ is $v=n\,2\,\cdots\,(n-2)\,(n-1)\,1$. It follows that $v$ is obtained from either $u_1$ or $u_2$ by a rotation of the integers in positions $1,n-1,n$, with possible sign changes of the integers in positions $1$ and $n-1$. The same happens between the last and first elements of the sequences generated in consecutive iterations of the {\bf for} cycle in Step 3 of the algorithm, $(\overline{i}\,\overline{n})$ and $(i\!+\!1\,n)$, respectively, which differ by the rotation of the letters in positions $i,i+1,n$, and sign changes of the letters in positions $i$ and $n$. The transitions  between Steps 5 and 6 are involutions that differ by two sign changes.
Finally, notice that property {\bf A1} is clearly satisfied by the sequence $\GCD(n)$.
\end{proof}

\begin{corollary}
The Gray code $\GCD(n)$ produced by Algorithm \ref{alg:typeD} is an Hamilton cycle in the restriction to involutions of the Cayley graph $G(\ss_n^D,T^D)$, with Hamming distance three, for $n\geq 3$.
\end{corollary}
\begin{proof}
We will show by induction on $n\geq 2$ that each element of $\GCD(n)$ is obtained multiplying on the right the previous element of the code by an element of $T^D$.
This is certainly true for $n=2$, where each element differ by the previous by two sign changes or a transposition. Fix $n\geq 3$ and assume the result is valid for integers less than $n$.
We start by noticing that by the proof of Corollary \ref{cor:hamtypeB}, each two consecutive elements of the code $\GCB(n-1)$ differ by a cycle of three letters in positions $i<j<k$  with the same sign, a transposition of two letters in positions $i<j$, with at most a sign change of the letter in position $i$, or one or two sign changes. This means that each involution generated in Step $2$ of Algorithm \ref{alg:typeD} is obtained by multiplying the previous element on the right by the generators $t_{j,k-1}\cdot t_{i,j-1}$, $t_{i,j}\cdot t_{j+1,n}\cdot t_{n-1,n}$, $t_{i,n}\cdot t_{j,n}$, or $t_{i,n}$. The same happens with each involution of $\GCD^F(n-2)\cdot (i\,n)$ and  $\overleftarrow{\GCD}^F(n-2)\cdot(\overline{i}\,\overline{n})$, generated in each iteration of the cycle {\bf for} of the algorithm. Finally, notice that by the proof of  Theorem
\ref{teor:algtypeD}, the involutions connecting each one of factors $\GCB^F(n-1)\cdot \widetilde{n}$, $\GCD^F(n-2)\cdot (i\,n)$ and $\overleftarrow{\GCD}^F(n-2)\cdot(\overline{i}\,\overline{n})$, for $i=1,\ldots, n-1$, are obtained by a rotation of three letters in positions $i<j<k$, with possibly sign changes of the letters in positions $i$ and $j$, or in positions $i$ and $k$. This is achieved by right multiplication of the correspondent involution by $t_{j,k-1}\cdot t_{i,j-1}$, or $t_{j,k-1}\cdot t_{i,j-1}\cdot t_{i,n}\cdot t_{j,n}$, or $t_{j,k-1}\cdot t_{i,j-1}\cdot t_{i,n}\cdot t_{k,n}$.
\end{proof}

\section{An optimal algorithm for the involutions of type $B$}\label{sec:optimal_algB}

In the previous section we have seen that the Gray Code presented for involutions of type $A$ is optimal with respect to Hamming distance. The same is not true for the codes presented for involutions of types $B$ and $D$. In \cite{grz1} it is proven that the minimal Hamming distance of a Gray Code for type $B$ involutions is two and a Gray Code with Hamming distance two is given. Is this section we will briefly describe that algorithm and, alongside with the description, we will construct the code for $\ss^B_4$ to better illustrate the steps of algorithm.
A modification of this algorithm would provide a Gray Code for  type $D$ involutions with Hamming distance two, which would be optimal.
We give such a construction for $\ss^D_4$ and $\ss^D_5$.
Before starting the description of the code, we introduce some notation for type $B$ involutions. Since $\ss^D_n\subset \ss^B_n$, the some notation is valid for both types.

Given $i\in[\pm n]$, we define the sign function $\sgn(i)=0$ if $i>0$, and $\sgn(i)=1$ otherwise.
A sign permutation $\pi$ can also be represented by the pair $(p,g)$, where
$p\in\mathcal{S}_n$ is defined by $p(i)=|\pi(i)|$ for all $i\in[n]$, and $g=g_1\cdots g_n\in B_n$, the set of binary words of length $n$, is defined by $g_i=\sgn(\pi(i))$ for all $i\in[n]$.
For example, $\pi=\overline{3}\,2\,\overline{1}\,6\,5\,4\,\overline{7}\in\mathcal{S}_7^B$ corresponds to the pair $(3216547,1010001)\in \mathcal{S}_7\times B_7$.

Since the elements of $\mathcal{S}_n^B$ are permutations of $[\pm n]$, we can write them in disjoint cycle form. In particular, if $\pi= (p,g)$ is an involution in $\mathcal{S}_n^B$, then $p$ is also an involution in $\mathcal{S}_n$, and
the cycle decomposition of $\pi= (p,g)$ is obtained by writing $p$ as the disjoint union
of transpositions
and position fixed points, {\it i.e.} integers $i$ for which $p(i)=i$, and then associating the respective signs.
For our running example $\pi= (3216547,1010001)$, we have $3216547=(13)(46)(2)(5)(7)\in\mathcal{S}_7$, and thus $\pi=(\overline{1}\,\overline{3})(46)(2)(5)(\overline{7})$.
When writing a transposition $(a\,b)$ we adopt the convention that $|a|<|b|$.
The integers $a$ and $b$ are called, respectively, the {\it opener} and {\it closure} of the transposition $(a\,b)$ and, as part of an involution, must have the same sign.
 The   common sign of both elements in a transposition is called a {\it paired sign}.

\vspace*{5mm}

Our approach is based on the successive application  of
the Binary Reflected Gray Code \cite{gray}, $BRGC_n$, on the set of all binary words that can be associated to a particular involution of $\ss_n$ to form a signed involution.
We recall now what is the Binary Reflected Gray Code. The $BRGC_1$ is the list $(0,1)$, and the $BRGC_{n+1}$ is obtainable by first listing $BRGC_n$, with each word
prefixed by $0$, and then listing the $BRGC_n$ in reverse order with each word prefixed by $1$.
For instance,
$BRGC_2=(00,01,11,10)$ and $BRGC_3=(000,001,011,010,$ $110,111,101,100).$
Each one of the $2^n$ elements of the sequence $BRGC_n$ differs from the previous one by only one bit, and the same is true for the last and the first elements of the sequence.

This second Gray code, designated by $OGCB_n$, for the involutions in $\ss^B_n$ is constructed  by layers.
Each layer consists in all involutions with a fix
 number of transpositions.
Given $n\in\mathbb{N}$, and $0\leq k\leq \lfloor n/2\rfloor$, let $L_{k}$ denote the set of all  involutions in $\ss_n$ having exactly $k$ transpositions. %, and let $$L_{\leq k}=\bigcup_{i=0}^kL_k$$ be the set of all involutions having at most $k$ transpositions.
For a sign involution $\sigma=(p,g)$, we say that $\sigma$ is in $L_k^B$ whenever $p\in L_k$. The only element in $L_{0}$ is the identity, and  any element in $L_0^B$ is written as $(id,g)$ with $g\in B_n$.

\begin{enumerate}
\item[($L_0^B$)] In the first stage we construct the code for the involutions in $L_0^B$ as the sequence $(id,g^1),\ldots, (id,g^{2^n})$, where $(g^1,\ldots, g^{2^n})=BRGC_n$.

\vspace*{5mm}

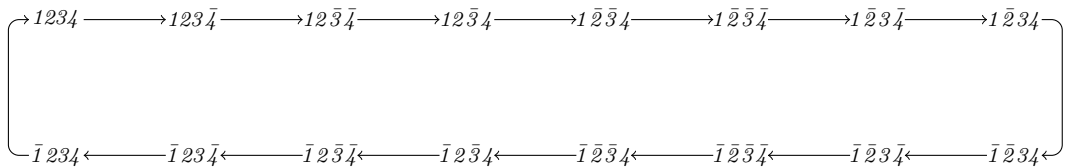
\begin{figure}[h]
  \begin{tikzpicture}[scale=0.9]

%%%%%%%%%%%
%%% L_0 %%%
%%%%%%%%%%%
%\node at (0,4.25) {$L_0^B$};

%%%\draw [rounded corners, dashed] (-0.25,1.5) rectangle (15,4.5);

%L_0(1234)
%\node at (7.5,3) {$\textbf{[1234]}$};
%
\draw [<-, rounded corners] (0.1,4) -- (-0.2,4) -- (-0.2, 2) -- (0.1, 2);
\node at (0.5,2) {{\tiny $\mathit{\bar{1}234}$}};
\draw [->, rounded corners] (2.1, 2) -- (.9, 2);
\node at (2.5,2) {{\tiny $\mathit{\bar{1}23\bar{4}}$}};
\draw [->, rounded corners] (4.1, 2) -- (2.9, 2);
\node at (4.5,2) {{\tiny $\mathit{\bar{1}2\bar{3}\bar{4}}$}};
\draw [->, rounded corners] (6.1, 2) -- (4.9, 2);
\node at (6.5,2) {{\tiny $\mathit{\bar{1}2\bar{3}4}$}};
\draw [->, rounded corners] (8.1, 2) -- (6.9, 2);
\node at (8.5,2) {{\tiny $\mathit{\bar{1}\bar{2}\bar{3}4}$}};
\draw [->, rounded corners] (10.1, 2) -- (8.9, 2);
\node at (10.5,2) {{\tiny $\mathit{\bar{1}\bar{2}\bar{3}\bar{4}}$}};
\draw [->, rounded corners] (12.1, 2) -- (10.9, 2);
\node at (12.5,2) {{\tiny $\mathit{\bar{1}\bar{2}3\bar{4}}$}};
\draw [->, rounded corners] (14.1, 2) -- (12.9, 2);
\node at (14.5,2) {{\tiny $\mathit{\bar{1}\bar{2}34}$}};

\node at (0.5,4) {{\tiny $\mathit{1234}$}};
\draw [->, rounded corners] (0.9, 4) -- (2.1, 4);
\node at (2.5,4) {{\tiny $\mathit{123\bar{4}}$}};
\draw [->, rounded corners] (2.9, 4) -- (4.1, 4);
\node at (4.5,4) {{\tiny $\mathit{12\bar{3}\bar{4}}$}};
\draw [->, rounded corners] (4.9, 4) -- (6.1, 4);
\node at (6.5,4) {{\tiny $\mathit{12\bar{3}4}$}};
\draw [->, rounded corners] (6.9, 4) -- (8.1, 4);
\node at (8.5,4) {{\tiny $\mathit{1\bar{2}\bar{3}4}$}};
\draw [->, rounded corners] (8.9, 4) -- (10.1, 4);
\node at (10.5,4) {{\tiny $\mathit{1\bar{2}\bar{3}\bar{4}}$}};
\draw [->, rounded corners] (10.9, 4) -- (12.1, 4);
\node at (12.5,4) {{\tiny $\mathit{1\bar{2}3\bar{4}}$}};
\draw [->, rounded corners] (12.9, 4) -- (14.1, 4);
\node at (14.5,4) {{\tiny $\mathit{1\bar{2}34}$}};
\draw [->, rounded corners] (14.9,4) -- (15.2,4) -- (15.2, 2) -- (14.9, 2);

  \end{tikzpicture}
  \caption{Gray code scheme for $L_0^B$ with $n=4$.}
  \label{fig:L_0, n=4}
\end{figure}

\vspace*{5mm}

\item[($L_1^B$)] After having the code for $L_{0}^B$, we construct the code for $L_{0}^B\cup L_{1}^B$ by inserting sequences of involutions in $L_{1}^B$ between two consecutive elements of $L_{0}^B$.  Each one of these sequences is of the form
    $$\left((p,h^1),\ldots,(p,h^{2^{n-1}})\right),$$
    where $p=(st)$ is a transposition in $L_1$ and $h=(h^1,\ldots,h^{2^{n-1}})$
    is a sequence of binary words of length $n$, obtained from the $BRGC_{n-1}$. The  position $\ell$ where the insertion of this sequence will occur  correspond to an involution
    $(id, g^{\ell})$ where the sign of the letters $s$ and $t$ is the same, and remain unchanged in the involution in position $\ell+1$.
    In \cite{grz2} a procedure was given to compute distinct places for each of these insertions, as well as the correspondent sequence of binary words.

For instance, with $n=4$ and $p=(12)$, the associated binary sequence  is $h=\big(1100,0000,0010,1110,1111,0011,0001,1101\big)$, which will be inserted in positions $\ell=9$  of $L_0^B$.
Figure \ref{figureB4_orelha1} shows the insertion of the involutions associated with transposition $(12)$ into $L_0^B$.

\vspace*{5mm}

\begin{figure}[H]
  \begin{tikzpicture}[scale=0.9]

%%%%%%%%%%%
%%% L_0 %%%
%%%%%%%%%%%
%\node at (0,4.25) {$L_0^B$};

%%%\draw [rounded corners, dashed] (-0.25,1.5) rectangle (15,4.5);

%L_0(1234)
%\node at (7.5,3) {$\textbf{[1234]}$};
%
\draw [<-, rounded corners] (0.1,4) -- (-0.2,4) -- (-0.2, 2) -- (0.1, 2);
\node at (0.5,2) {{\tiny $\mathit{\bar{1}234}$}};
\draw [->, rounded corners] (2.1, 2) -- (.9, 2);
\node at (2.5,2) {{\tiny $\mathit{\bar{1}23\bar{4}}$}};
\draw [->, rounded corners] (4.1, 2) -- (2.9, 2);
\node at (4.5,2) {{\tiny $\mathit{\bar{1}2\bar{3}\bar{4}}$}};
\draw [->, rounded corners] (6.1, 2) -- (4.9, 2);
\node at (6.5,2) {{\tiny $\mathit{\bar{1}2\bar{3}4}$}};
\draw [->, rounded corners] (8.1, 2) -- (6.9, 2);
\node at (8.5,2) {{\tiny $\mathit{\bar{1}\bar{2}\bar{3}4}$}};
\draw [->, rounded corners] (10.1, 2) -- (8.9, 2);
\node at (10.5,2) {{\tiny $\mathit{\bar{1}\bar{2}\bar{3}\bar{4}}$}};
\draw [->, rounded corners] (12.1, 2) -- (10.9, 2);
\node at (12.5,2) {{\tiny $\mathit{\bar{1}\bar{2}3\bar{4}}$}};
%\draw [->, rounded corners] (14.1, 2) -- (12.9, 2);
\node at (14.5,2) {{\tiny $\mathit{\bar{1}\bar{2}34}$}};

\node at (0.5,4) {{\tiny $\mathit{1234}$}};
\draw [->, rounded corners] (0.9, 4) -- (2.1, 4);
\node at (2.5,4) {{\tiny $\mathit{123\bar{4}}$}};
\draw [->, rounded corners] (2.9, 4) -- (4.1, 4);
\node at (4.5,4) {{\tiny $\mathit{12\bar{3}\bar{4}}$}};
\draw [->, rounded corners] (4.9, 4) -- (6.1, 4);
\node at (6.5,4) {{\tiny $\mathit{12\bar{3}4}$}};
\draw [->, rounded corners] (6.9, 4) -- (8.1, 4);
\node at (8.5,4) {{\tiny $\mathit{1\bar{2}\bar{3}4}$}};
\draw [->, rounded corners] (8.9, 4) -- (10.1, 4);
\node at (10.5,4) {{\tiny $\mathit{1\bar{2}\bar{3}\bar{4}}$}};
\draw [->, rounded corners] (10.9, 4) -- (12.1, 4);
\node at (12.5,4) {{\tiny $\mathit{1\bar{2}3\bar{4}}$}};
\draw [->, rounded corners] (12.9, 4) -- (14.1, 4);
\node at (14.5,4) {{\tiny $\mathit{1\bar{2}34}$}};
\draw [->, rounded corners] (14.9,4) -- (15.2,4) -- (15.2, 2) -- (14.9, 2);

%%%%%%%%%%%%L_1

\node at (12,0) {\tiny $21\bar{3}\bar{4}$};
\node at (13,0) {\tiny $\bar{2}\bar{1}\bar{3}\bar{4}$};
\node at (14,0) {\tiny $\bar{2}\bar{1}\bar{3}4$};
\node at (15,0) {\tiny $21\bar{3}4$};
\node at (12,1) {\tiny $213\bar{4}$};
\node at (13,1) {\tiny $\bar{2}\bar{1}3\bar{4}$};
\node at (14,1) {\tiny $\bar{2}\bar{1}34$};
\node at (15,1) {\tiny $2134$};
\draw [rounded corners,gray,->]  (13.2,1.2) -- (13.2,2) -- (12.9,2);
\draw [rounded corners,gray,->]  (14.1,2) -- (13.9,2) -- (13.9,1.2);
\draw [->] (14.4,1) -- (14.6,1);
\draw [->] (15,0.8) -- (15,0.2);
\draw [->] (14.6,0) -- (14.4,0);
\draw [->] (13.6,0) -- (13.4,0);
\draw [->] (12,.2) -- (12,.8);
\draw [->] (12.4,1) -- (12.6,1);
\draw [->] (12.6,0) -- (12.4,0);

  \end{tikzpicture}
  \caption{Insertion of the sign involutions associated with $(12)$ into $L_0^B$.}
  \label{figureB4_orelha1}
\end{figure}
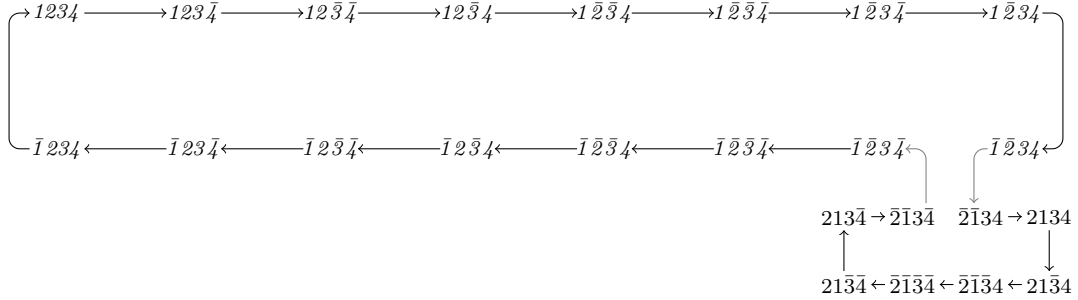

\vspace*{5mm}

\item[($L_2^B$)] Next, we have to insert the elements of $L_2^B$ into the code for $L_0^B\cup L_1^B$ computed in the previous step. Each permutation in $L_2$ is the  product of two transpositions, which will be sorted by the lexicographic order of the openers. As in the previous case, for each  permutation $p'\in L_2$, we can compute a sequence of binary words $(f^1,\ldots,f^{2^{n-2}})$ such that
    $\left((p',f^1),\ldots,(p',f^{2^{n-2}})\right)$ contains all sign involutions associated with $p$, as well as distinct places places in the code for $L_0^B\cup L_1^B$ to insert such sequences, satisfying the same property as in the previous layer.

In our running example, with $n=4$, $p=(12)$ and $p'=p\cdot (34)$,  the associated binary code is $(1111,1100,0000,0011)$.
The sequence $$\big((p',1111),(p',1100), (p',0000), (p',0011)\big)$$ will be inserted in position $\ell=5$, that is between involutions $(p,h^5)$ and $(p,h^6)$.
Figure \ref{figureB4_orelha2} shows the insertion of the  sign involutions associated with  $p\cdot(34)$ into $L_0^B\cup L_1^B$.

\vspace*{5mm}

\begin{figure}[H]
  \begin{tikzpicture}[scale=1]

%%%%%%%%%%%
%%% L_0 %%%
%%%%%%%%%%%
%\node at (0,4.25) {$L_0^B$};

%%%\draw [rounded corners, dashed] (-0.25,1.5) rectangle (15,4.5);

%L_0(1234)
%\node at (7.5,3) {$\textbf{[1234]}$};
%
\node at (12.5,2) {{\tiny $\mathit{\bar{1}\bar{2}3\bar{4}}$}};
%\draw [->, rounded corners] (14.1, 2) -- (12.9, 2);
\node at (14.5,2) {{\tiny $\mathit{\bar{1}\bar{2}34}$}};
\draw [->, rounded corners] (12, 2) -- (11.4,2);
\node at (11,2) {$\cdots$};
\draw [->, rounded corners] (15.6, 2) -- (15,2);
\node at (16.05,2) {$\cdots$};

%%%%%%%%%%%%L_1

\node at (12,0) {\tiny $21\bar{3}\bar{4}$};
\node at (13,0) {\tiny $\bar{2}\bar{1}\bar{3}\bar{4}$};
\node at (14,0) {\tiny $\bar{2}\bar{1}\bar{3}4$};
\node at (15,0) {\tiny $21\bar{3}4$};
\node at (12,1) {\tiny $213\bar{4}$};
\node at (13,1) {\tiny $\bar{2}\bar{1}3\bar{4}$};
\node at (14,1) {\tiny $\bar{2}\bar{1}34$};
\node at (15,1) {\tiny $2134$};
\draw [rounded corners,gray,->]  (13.2,1.2) -- (13.2,2) -- (12.9,2);
\draw [rounded corners,gray,->]  (14.1,2) -- (13.9,2) -- (13.9,1.2);
\draw [rounded corners,->] (14.4,1) -- (14.6,1);
\draw [rounded corners,->] (15,0.8) -- (15,0.2);
\draw [rounded corners,->] (14.6,0) -- (14.4,0);
\draw [rounded corners,->] (13.6,0) -- (13.4,0);
\draw [rounded corners,->] (12,.2) -- (12,.8);
\draw [rounded corners,->] (12.4,1) -- (12.6,1);

%%%%%%%%%%%%%%%L_2

\node at (12,-1) {\tiny $21\bar{4}\bar{3}$};
\node at (13,-1) {\tiny $\bar{2}\bar{1}\bar{4}\bar{3}$};
\node at (13,-1.5) {\tiny $\bar{2}\bar{1}43$};
\node at (12,-1.5) {\tiny $2143$};
%\node at (12.5,-1.25) {$2\textbf{1}4\textbf{3}$};
%\draw (12,-1.5) to (13,-1.5);
\draw [rounded corners,->] (13.3,-1) -- (13.5,-1) -- (13.5,-1.5) -- (13.3,-1.5);
\draw [rounded corners,->] (11.7,-1.5) -- (11.5,-1.5) -- (11.5,-1) -- (11.7,-1);
\draw [rounded corners,->] (12.6,-1.5) -- (12.4,-1.5);
\draw [rounded corners,->,gray] (12,-.8) -- (12,-.2);
\draw [rounded corners,->,gray] (13,-.2) -- (13,-.8);

  \end{tikzpicture}
  \caption{Insertion of the sign involutions associated with $(12)(34)$ into $L_0^B\cup L_1^B$.}
  \label{figureB4_orelha2}
\end{figure}
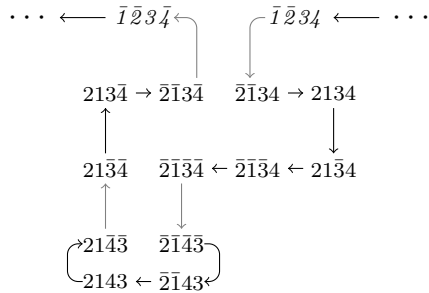

\vspace*{5mm}

\item[($L_{k}^B$)] The insertion of the involutions of $L_k^B$ into the code for $L_0^B\cup\cdots\cup L_{k-1}^B$ is similar to the insertion of the  involutions of $L_2^B$ into the previous layer.  Each sequence is associated to an element of
    $$L_{k}(q):=\{q\cdot (ij)\,:\, m(q)<i<j,\, q_i=i,\,q_j=j\}\subseteq L_{k},$$
 where $m(q)$ is the largest opener amongst all transpositions in the cycle decomposition of $q$, and  $m(q)=0$ when $q$ is the identity. Each sequence   will be inserted into the correspondence sequence associated with the involution $q$, following the same rules as before. Notice that the sets $L_{k}(q)$ form a partition of $L_k$.
\end{enumerate}

\medskip

The code $OGCB_n$ is obtained when the last layer is inserted.

\begin{theorem}\label{teo:hamtypeB2}
The code $OGCB_n$  produces an Hamilton cycle in the restriction to involutions of the Cayley graph $G(\ss_n^B,X_1^B\cup X_2^B)$, with Hamming distance two, for $n\geq 3$.
This is the minimal Hamming distance of any Gray code for $I_n^B$.
\end{theorem}

\begin{proof}
In \cite{grz1} it was proven that the sequence produced by the code $OGCB_n$ is a cyclic Gray code with Hamming distance two for $\ss_n^B$, where two consecutive involutions differ by one sign change, two sign changes - in this case, a paired sign -, or a transposition without any sign change.
This means that each involution, in the code, is obtained by multiplying the previous element on the right by one reflections of the set $X_1^B\cup X_2^B$.
Note that the Hamming distance for any code can only be one if there are only single sign changes between consecutive involutions. Thus this code has the minimal Hamming distance for $I_n^B$.
\end{proof}

Table \ref{table:graycodesign4dist2} shows the Gray code with Hamming distance 2 for the sign involutions in $I_4^B$ using the algorithm defined above. The code  should be read down columns, from left to right.

\begin{table}[h]
 \centering
\begin{tabular}{|c|c|c|c|c|c|c|c|}
\hline
$1\,2\,3\,4$ & $\bar{1}\,3\,2\,\bar{4}$ & $1\,4\,3\,2$ & $\bar{2}\,\bar{1}\,\bar{4}\,\bar{3}$ & $4\,2\,3\,1$ & $\bar{1}\,\bar{2}\,\bar{3}\,\bar{4}$ & $3\,2\,1\,4$ & $\bar{1}\,\bar{2}\,\bar{4}\,\bar{3}$\\
$1\,2\,3\,\bar{4}$& $1\,3\,2\,\bar{4}$ & $1\,\bar{4}\,3\,\bar{2}$ & $\bar{2}\,\bar{1}\,4\,3$ & $\bar{4}\,2\,3\,\bar{1}$ & $\bar{3}\,\bar{2}\,\bar{1}\,\bar{4}$& $3\,\bar{2}\,1\,4$ &$\bar{1}\,\bar{2}\,4\,3$\\
$1\,2\,\bar{3}\,\bar{4}$&$1\,\bar{3}\,\bar{2}\,\bar{4}$& $1\,\bar{2}\,3\,\bar{4}$ & $2\,1\,4\,3$ & $\bar{4}\,2\,\bar{3}\,\bar{1}$ & $3\,\bar{2}\,1\,\bar{4}$& $\bar{3}\,\bar{2}\,\bar{1}\,4$ &$1\,\bar{2}\,4\,3$\\
$1\,2\,\bar{3}\,4$&$1\,\bar{2}\,\bar{3}\,\bar{4}$& $1\,\bar{2}\,3\,4$ & $2\,1\,\bar{4}\,\bar{3}$ & $4\,2\,\bar{3}\,1$& $3\,2\,1\,\bar{4}$& $\bar{1}\,\bar{2}\,\bar{3}\,4$ &$1\,\bar{2}\,\bar{4}\,\bar{3}$\\
$1\,\bar{2}\,\bar{3}\,4$&$1\,\bar{4}\,\bar{3}\,\bar{2}$& $\bar{1}\,\bar{2}\,3\,4$ & $2\,1\,\bar{3}\,\bar{4}$  & $4\,\bar{2}\,\bar{3}\,1$&$\bar{3}\,2\,\bar{1}\,\bar{4}$& $\bar{1}\,2\,\bar{3}\,4$&$1\,2\,\bar{4}\,\bar{3}$\\
$1\,\bar{3}\,\bar{2}\,4$&$1\,4\,\bar{3}\,2$& $\bar{2}\,\bar{1}\,3\,4$ & $2\,1\,3\,\bar{4}$ & $4\,\bar{3}\,\bar{2}\,1$& $\bar{3}\,2\,\bar{1}\,4$& $\bar{1}\,2\,\bar{3}\,\bar{4}$&$1\,2\,4\,3$\\
$1\,3\,2\,4$&$\bar{1}\,4\,\bar{3}\,2$& $2\,1\,3\,4$ & $\bar{2}\,\bar{1}\,3\,\bar{4}$ & $4\,3\,2\,1$ & $\bar{3}\,4\,\bar{1}\,2$& $\bar{1}\,2\,3\,\bar{4}$ &\\
$\bar{1}\,3\,2\,4$& $\bar{1}\,\bar{4}\,\bar{3}\,\bar{2}$ & $2\,1\,\bar{3}\,4$ &$\bar{1}\,\bar{2}\,3\,\bar{4}$ & $\bar{4}\,3\,2\,\bar{1}$ & $\bar{3}\,\bar{4}\,\bar{1}\,\bar{2}$& $\bar{1}\,2\,3\,4$ &\\
$\bar{1}\,\bar{3}\,\bar{2}\,4$& $\bar{1}\,\bar{4}\,3\,\bar{2}$& $\bar{2}\,\bar{1}\,\bar{3}\,4$ & $\bar{4}\,\bar{2}\,3\,\bar{1}$ & $\bar{4}\,\bar{3}\,\bar{2}\,\bar{1}$ & $3\,\bar{4}\,1\,\bar{2}$& $\bar{1}\,2\,4\,3$ &\\
$\bar{1}\,\bar{3}\,\bar{2}\,\bar{4}$& $\bar{1}\,4\,3\,2$& $\bar{2}\,\bar{1}\,\bar{3}\,\bar{4}$ & $4\,\bar{2}\,3\,1$ & $\bar{4}\,\bar{2}\,\bar{3}\,\bar{1}$& $3\,4\,1\,2$ & $\bar{1}\,2\,\bar{4}\,\bar{3}$&\\
\hline
\end{tabular}
\caption{Gray code for $I_4^B$ with Hamming distance 2.}\label{table:graycodesign4dist2}
\end{table}

\subsection{Towards an optimal algorithm for the involutions of type $D$}

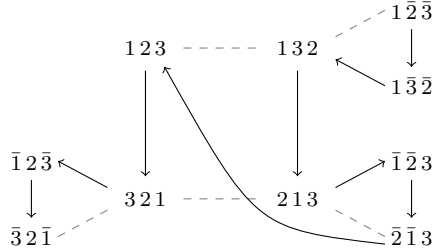
\begin{figure}[H]
  \begin{tikzpicture}[scale=1]

\node at (-1,1) {\tiny{$1\,2\,3$}};
\node at (1,1) {\tiny{$1\,3\,2$}};
\node at (1,-1) {\tiny{$2\,1\,3$}};
\node at (-1,-1) {\tiny{$3\,2\,1$}};
\draw[gray,dashed]  (-0.5, 1) -- (0.5, 1);
\draw[gray,dashed]  (-0.5, -1) -- (0.5, -1);
\draw[->]  (-1, 0.7) -- (-1, -0.7);
\draw[->]  (1, 0.7) -- (1, -0.7);

\node at (2.5,1.5) {\tiny{$1\,\bar{2}\,\bar{3}$}};
\node at (2.5,0.5) {\tiny{$1\,\bar{3}\,\bar{2}$}};
\draw[gray,dashed]  (1.5, 1.15) -- (2.15, 1.5);
\draw[<-]  (1.5, 0.85) -- (2.15, 0.5);
\draw[->]  (2.5, 1.25) -- (2.5, 0.75);

\node at (2.5,-1.5) {\tiny{$\bar{2}\,\bar{1}\,3$}};
\node at (2.5,-0.5) {\tiny{$\bar{1}\,\bar{2}\,3$}};
\draw[gray,dashed]  (1.5, -1.15) -- (2.15, -1.5);
\draw[->]  (1.5, -0.85) -- (2.15, -0.5);
\draw[<-]  (2.5, -1.25) -- (2.5, -0.75);

\node at (-2.5,-1.5) {\tiny{$\bar{3}\,2\,\bar{1}$}};
\node at (-2.5,-0.5) {\tiny{$\bar{1}\,2\,\bar{3}$}};
\draw[gray,dashed]  (-1.5, -1.15) -- (-2.15, -1.5);
\draw[->]  (-1.5, -0.85) -- (-2.15, -0.5);
\draw[<-]  (-2.5, -1.25) -- (-2.5, -0.75);

\draw[->]  (2.15, -1.6) .. controls (.45,-1.4) ..  (-0.75, 0.75);
  \end{tikzpicture}
  \caption{An Hamilton path in $G(I_3^D)$ with Hamming distance 2.}
  \label{figure:graphB3}
\end{figure}

In section \ref{section:typeD} we constructed an Hamilton cycle with Hamming distance 3 for the involutions in $I_n^D$. In this section we analyze  the existence of such Hamilton cycle with distance 2, which would be optimal since is clear that such code cannot have distance 1. So, consider the graph $G(I_n^D)$, where two involutions are connected if their Hamming distance is 2. When $n=2$ the Hamming distance between any two involutions in  $I_2^D$ is 2, and thus $G(I_2^D)$ has an Hamilton cycle with Hamming distance 2. The graph $G(I_3^D)$ is displayed  in Figure \ref{figure:graphB3}, where we have omitted all but three of the edges  linking the identity $123$ to all other involutions. It is easy to check that there is no Hamilton cycle in $G(I_3^D)$, but there are Hamilton paths. One Hamilton path  starts with the involution $1\,\bar{2}\,\bar{3}$, and  can be read off of the graph  by following the solid edges:
$$\left(1\,\bar{2}\,\bar{3},\, 1\,\bar{3}\,\bar{2},\, 1\,3\,2,\,2\,1\,3,\,\bar{1}\,\bar{2}\,3,\,
\bar{2}\,\bar{1}\,3,\,1\,2\,3,\, 3\,2\,1,\,\bar{1}\,2\,\bar{3},\,\bar{3}\,2\,\bar{1}\right).$$

\vspace*{5mm}

Our approach for the general case is based on the successive application  of
 a Binary  Code, $BCE_n$, on the set of all $2^{n-1}$ binary words with an even number of $1's$ that can be associated to a particular involution of $\ss_n$ to form an involution of type $D$. This Binary  Code is constructed by recursion for $n\geq 2$, triggered by $BCE_2=(00,\, 11)$, as follows: if
 $$BCE_{n-1}=(u_1,u_2,\ldots,u_k)\quad\text{ and }\quad BRGC_{n-2}=(v_1,v_2,\ldots,v_{k}),$$
 where $k=2^{n-2}$,
 then let
 $$BCE_n=(0\cdot u_1,0\cdot u_2,\ldots, 0\cdot u_k, 10\cdot v_k, 11\cdot v_{k-1}, 10\cdot v_{k-2}, 11\cdot v_{k-3},\ldots, 11\cdot v_1).$$

 Each one of the $2^{n-1}$ elements of the sequence $BCE_n$ differs from the previous one by two bits, and the same is true for the last and the first elements of the sequence.
For instance, the $BCE_n$ code for $n=3,4,5$ is given in Table \ref{table:BCE345}.

\begin{table}[h]
 \centering
\begin{tabular}{|c|c|c|}
\hline
$BCE_3$&$BCE_4$&$BCE_5$\\
\hline
000&0000&00000\\
011&0011&00011\\
101&0101&00101\\
110&0110&00110\\
&1010&01010\\
&1111&01111\\
&1001&01001\\
&1100&01100\\
&&10100\\
&&11101\\
&&10111\\
&&11110\\
&&10010\\
&&11011\\
&&10001\\
&&11000\\
\hline
\end{tabular}
\caption{$BCE_n$, for $n=3,4,5$.}
  \label{table:BCE345}
\end{table}

The construction of an Hamilton cycle for $G(I_4^D)$ with distance 2 can be obtained as follows.
We start by constructing a path with distance 2 for the sign involution in the set $\{(p,g): g\in B_4\}\cap I_4^D$, for each $p\in L_2$:
\begin{align*}
L_2^{(12)(34)}&=(2\,1\,\bar{4}\,\bar{3},\,\bar{2}\,\bar{1}\,\bar{3}\,\,\bar{4},\bar{2}\,\bar{1}\,4\,3,\,2\,1\,4\,3),\\
L_2^{(13)(24)}&=(3\,4\,1\,2,\,\bar{3}\,4\,\bar{1}\,2,\,\bar{3}\,\bar{4}\,\bar{1}\,\bar{2},\,3\,\bar{4}\,1\,\bar{2}),\\
L_2^{(14)(23)}&=(4\,\bar{3}\,\bar{2}\,1,\,4\,3\,2\,1,\,\bar{4}\,3\,2\,\bar{1},\,\bar{4}\,\bar{3}\,\bar{2}\,\bar{1}).\\
\end{align*}

Each one of these sequences will be inserted between the paths with distance 2 for two sets
$\{(p,g):  g\in B_4\}\cap I_4^D$, with $p\in L_1$, such that the elements linking the three paths differ by a single transposition of two letters:

\begin{align*}
L_1^{(12)}\cdot L_2^{(12)(34)}\cdot L_1^{(34)}&=(2\,1\,3\,4,\,\bar{2}\,\bar{1}\,3\,4,\,\bar{2}\,\bar{1}\,\bar{3}\,\bar{4},\,2\,1\,\bar{3}\,\bar{4},\,L_2^{(12)(34)},\,1\,2\,4\,3,\,\bar{1}\,\bar{2}\,4\,3,\,\bar{1}\,\bar{2}\,\bar{4}\,\bar{3},\,1\,2\,\bar{4}\,\bar{3}),\\
L_1^{(13)}\cdot L_2^{(13)(24)}\cdot L_1^{(24)}&=(\bar{3}\,2\,\bar{1}\,4,\,\bar{3}\,\bar{2}\,\bar{1}\,\bar{4},\,3\,\bar{2}\,1\,\bar{4},\,3\,2\,1\,4,\,L_2^{(13)(24)},\,1\,\bar{4}\,3\,\bar{2},\,1\,4\,3\,2,\,\bar{1}\,4\,\bar{3}\,2,\,\bar{1}\,\bar{4}\,\bar{3}\,\bar{2}),\\
L_1^{(23)}\cdot L_2^{(23)(14)}\cdot L_1^{(14)}&=(\bar{1}\,\bar{3}\,\bar{2}\,\bar{4},\,\bar{1}\,3\,2\,\bar{4},\,1\,3\,2\,4,\,1\,\bar{3}\,\bar{2}\,4,L_2^{(14)(23)},\,\bar{4}\,\bar{2}\,\bar{3}\,\bar{1},\,4\,\bar{2}\,\bar{3}\,1,\,4\,2\,3\,1,\,\bar{4}\,2\,3\,\bar{1}).\\
\end{align*}

The final step consists in the insertion of these paths into the level $L_0$ path $$((1234,g_1),\ldots, (1234,g_8)),$$ where $(g_1,\ldots,g_8)=BCE_4$, such that the elements linking the different paths differ by a single transposition of two letters. The resulting code is an Hamilton cycle with distance 2 in the graph $G(I_4^D)$, where two consecutive involutions differ either by two sign changes or by a transposition of two letters:
\begin{align*}
I_4^D=&\big(1\,2\,3\,4,\,L_1^{(12)}\cdot L_2^{(12)(34)}\cdot L_1^{(34)},\,1\,2\,\bar{3}\,\bar{4},\,1\,\bar{2}\,3\,\bar{4},\,1\,\bar{2}\,\bar{3}\,4,\,\bar{1}\,2\,\bar{3}\,4,\\
&L_1^{(13)}\cdot L_2^{(13)(24)}\cdot L_1^{(24)},\,\bar{1}\,\bar{2}\,\bar{3}\,\bar{4},\,L_1^{(23)}\cdot L_2^{(23)(14)}\cdot L_1^{(14)},\,\bar{1}\,2\,3\,\bar{4},\,\bar{1}\,\bar{2}\,3\,4\big).
\end{align*}

The same process can be used to obtain an Hamilton cycle for $G(I_5^D)$ with distance 2. Start with  paths the $L_i^p$ with distance 2 for the sign involutions in the sets
$$\{(p,g): g\in BRGC_5\}\cap I_5^D,$$ for each $p\in L_i$, $i=2,1$, displayed in Table \ref{table:GI5}, and
$$L_0=\{(1234,g_1),\ldots,(1234,g_{16}):(g_1,\ldots,g_{16})=BCE_5\}.$$
Then, insert each path $L_2^p$ into a path $L_1^{p'}$, and finally each one of these paths into $L_0$, such that the words linking the paths differ by a single transposition of two letters, as indicated in the Table \ref{table:GI5}, which should be read top to bottom, starting in the left column.
\begin{table}[h]
 \centering
\begin{tabular}{|c|l|c|l|c|l|c|l|}
\hline
word&layer&word&layer&word&layer&word&layer\\
\hline
$1\,2\,3\,4\,5$ & $L_0$ & $1\,2\,\bar{4}\,\bar{3}\,5$ & $L_1^{(3 4)}$ & $\bar{1}\,2\,\bar{3}\,\bar{4}\,\bar{5}$ & $L_0$ & $\bar{5}\,2\,4\,3\,\bar{1}$ & $L_2^{(15) (34)}$\\
$1\,2\,3\,5\,4$ & $L_1^{(4 5)}$ & $1\,5\,\bar{4}\,\bar{3}\,2$ & $L_2^{(3 4)(25)}$ & $\bar{1}\,\bar{2}\,\bar{3}\,\bar{4}\,5$ & $L_0$ & $\bar{5}\,2\,\bar{4}\,\bar{3}\,\bar{1}$ & $L_2^{(15) (34)}$\\
$1\,2\,3\,\bar{5}\,\bar{4}$ & $L_1^{(4 5)}$ & $1\,5\,4\,3\,2$ & $L_2^{(3 4)(25)}$ & $\bar{4}\,\bar{2}\,\bar{3}\,\bar{1}\,5$ & $L_1^{(1 4)}$ & $\bar{5}\,2\,\bar{3}\,\bar{4}\,\bar{1}$ & $L_1^{(15)}$\\
$1\,\bar{2}\,\bar{3}\,\bar{5}\,\bar{4}$ & $L_1^{(4 5)}$ & $1\,\bar{5}\,4\,3\,\bar{2}$ & $L_2^{(3 4)(25)}$ & $\bar{4}\,\bar{3}\,\bar{2}\,\bar{1}\,5$ & $L_2^{(1 4) (2 3)}$ & $\bar{5}\,\bar{2}\,\bar{3}\,4\,\bar{1}$ & $L_1^{(15)}$\\
$1\,\bar{2}\,\bar{3}\,5\,4$ & $L_1^{(4 5)}$ & $1\,\bar{5}\,\bar{4}\,\bar{3}\,\bar{2}$ & $L_2^{(3 4)(25)}$ & $\bar{4}\,3\,2\,\bar{1}\,5$ & $L_2^{(1 4) (2 3)}$ & $5\,\bar{2}\,\bar{3}\,4\,1$ & $L_1^{(15)}$\\
$1\,\bar{3}\,\bar{2}\,5\,4$ & $L_2^{(4 5) (2 3)}$ & $1\,\bar{5}\,\bar{3}\,\bar{4}\,\bar{2}$ & $L_1^{(2 5)}$ & $4\,3\,2\,1\,5$ & $L_2^{(1 4) (2 3)}$ & $5\,\bar{3}\,\bar{2}\,4\,1$ & $L_2^{(15) (23)}$\\
$1\,\bar{3}\,\bar{2}\,\bar{5}\,\bar{4}$ & $L_2^{(4 5) (2 3)}$ & $1\,5\,\bar{3}\,\bar{4}\,2$ & $L_1^{(2 5)}$ & $4\,\bar{3}\,\bar{2}\,1\,5$ & $L_2^{(1 4) (2 3)}$ & $\bar{5}\,\bar{3}\,\bar{2}\,4\,\bar{1}$ & $L_2^{(15) (23)}$\\
$1\,3\,2\,\bar{5}\,\bar{4}$ & $L_2^{(4 5) (2 3)}$ & $1\,5\,3\,4\,2$ & $L_1^{(2 5)}$ & $4\,\bar{2}\,\bar{3}\,1\,5$ & $L_1^{(14)}$ & $\bar{5}\,3\,2\,4\,\bar{1}$ & $L_2^{(15) (23)}$\\
$1\,3\,2\,5\,4$ & $L_2^{(4 5) (2 3)}$ & $1\,\bar{5}\,3\,4\,\bar{2}$ & $L_1^{(2 5)}$ & $4\,2\,\bar{3}\,1\,\bar{5}$ & $L_1^{(14)}$ & $5\,3\,2\,4\,1$ & $L_2^{(15) (23)}$\\
$1\,3\,2\,4\,5$ & $L_1^{(2 3)}$ & $1\,\bar{2}\,3\,4\,\bar{5}$ & $L_0$ & $4\,2\,\bar{5}\,1\,\bar{3}$ & $L_2^{(1 4) (3 5)}$ & $5\,2\,3\,4\,1$ & $L_1^{(15)}$\\
$1\,\bar{3}\,\bar{2}\,4\,5$ & $L_1^{(23)}$ & $1\,\bar{2}\,\bar{3}\,4\,5$ & $L_0$ & $4\,2\,5\,1\,3$ & $L_2^{(1 4) (3 5)}$ & $\bar{5}\,2\,3\,4\,\bar{1}$ & $L_1^{(15)}$\\
$1\,\bar{3}\,\bar{2}\,\bar{4}\,\bar{5}$ & $L_1^{(23)}$ & $\bar{1}\,2\,\bar{3}\,4\,5$ & $L_0$ & $\bar{4}\,2\,5\,\bar{1}\,3$ & $L_2^{(1 4) (3 5)}$ & $\bar{1}\,2\,3\,4\,\bar{5}$ & $L_0$\\
$1\,3\,2\,\bar{4}\,\bar{5}$ & $L_1^{(23)}$ & $\bar{3}\,2\,\bar{1}\,4\,5$ & $L_1^{(1 3)}$ & $\bar{4}\,2\,\bar{5}\,\bar{1}\,\bar{3}$ & $L_2^{(1 4) (3 5)}$ & $\bar{1}\,\bar{2}\,3\,4\,5$ & $L_0$\\
$1\,2\,3\,\bar{4}\,\bar{5}$ & $L_0$ & $\bar{3}\,5\,\bar{1}\,4\,2$ & $L_2^{(1 3) (2 5)}$ & $\bar{4}\,2\,\bar{3}\,\bar{1}\,\bar{5}$ & $L_1^{(14)}$ & $\bar{2}\,\bar{1}\,3\,4\,5$ & $L_1^{(1 2)}$\\
$1\,2\,\bar{3}\,4\,\bar{5}$ & $L_0$ & $\bar{3}\,\bar{5}\,\bar{1}\,4\,\bar{2}$ & $L_2^{(1 3) (2 5)}$ & $\bar{4}\,\bar{2}\,3\,\bar{1}\,\bar{5}$ & $L_1^{(14)}$ & $\bar{2}\,\bar{1}\,4\,3\,5$ & $L_2^{(1 2) (3 4)}$\\
$1\,2\,\bar{3}\,\bar{4}\,5$ & $L_0$ & $3\,\bar{5}\,1\,4\,\bar{2}$ & $L_2^{(1 3) (2 5)}$ & $4\,\bar{2}\,3\,1\,\bar{5}$ & $L_1^{(14)}$ & $2\,1\,4\,3\,5$ & $L_2^{(1 2) (3 4)}$\\
$1\,\bar{2}\,3\,\bar{4}\,5$ & $L_0$ & $3\,5\,1\,4\,2$ & $L_2^{(1 3) (2 5)}$ & $4\,\bar{5}\,3\,1\,\bar{2}$ & $L_2^{(1 4) (25)}$ & $2\,1\,\bar{4}\,\bar{3}\,5$ & $L_2^{(1 2) (3 4)}$\\
$1\,\bar{2}\,5\,\bar{4}\,3$ & $L_1^{(3 5)}$ & $3\,2\,1\,4\,5$ & $L_1^{(13)}$ & $\bar{4}\,\bar{5}\,3\,\bar{1}\,\bar{2}$ & $L_2^{(1 4) (25)}$ & $\bar{2}\,\bar{1}\,\bar{4}\,\bar{3}\,5$ & $L_2^{(1 2) (3 4)}$\\
$1\,\bar{2}\,\bar{5}\,\bar{4}\,\bar{3}$ & $L_1^{(3 5)}$ & $3\,4\,1\,2\,5$ & $L_2^{(13) (24)}$ & $\bar{4}\,5\,3\,\bar{1}\,2$ & $L_2^{(1 4) (25)}$ & $\bar{2}\,\bar{1}\,\bar{3}\,\bar{4}\,5$ & $L_1^{(12)}$\\
$1\,2\,\bar{5}\,4\,\bar{3}$ & $L_1^{(3 5)}$ & $\bar{3}\,4\,\bar{1}\,2\,5$ & $L_2^{(13) (24)}$ & $4\,5\,3\,1\,2$ & $L_2^{(1 4) (25)}$ & $2\,1\,\bar{3}\,\bar{4}\,5$ & $L_1^{(12)}$\\
$1\,2\,5\,4\,3$ & $L_1^{(3 5)}$ & $\bar{3}\,\bar{4}\,\bar{1}\,\bar{2}\,5$ & $L_2^{(13) (24)}$ & $4\,2\,3\,1\,5$ & $L_1^{(14)}$ & $2\,1\,3\,\bar{4}\,\bar{5}$ & $L_1^{(12)}$\\
$1\,4\,5\,2\,3$ & $L_2^{(3 5)(24)}$ & $3\,\bar{4}\,1\,\bar{2}\,5$ & $L_2^{(13) (24)}$ & $\bar{4}\,2\,3\,\bar{1}\,5$ & $L_1^{(14)}$ & $2\,1\,3\,\bar{5}\,\bar{4}$ & $L_2^{(1 2) (45)}$\\
$1\,4\,\bar{5}\,2\,\bar{3}$ & $L_2^{(3 5)(24)}$ & $3\,\bar{2}\,1\,\bar{4}\,5$ & $L_1^{(13)}$ & $\bar{1}\,2\,3\,\bar{4}\,5$ & $L_0$ & $2\,1\,3\,5\,4$ & $L_2^{(1 2) (45)}$\\
$1\,\bar{4}\,\bar{5}\,\bar{2}\,\bar{3}$ & $L_2^{(3 5)(24)}$ & $\bar{3}\,\bar{2}\,\bar{1}\,\bar{4}\,5$ & $L_1^{(13)}$ & $\bar{1}\,\bar{2}\,3\,\bar{4}\,\bar{5}$ & $L_0$ & $\bar{2}\,\bar{1}\,3\,5\,4$ & $L_2^{(1 2) (45)}$\\
$1\,\bar{4}\,5\,\bar{2}\,3$ & $L_2^{(3 5)(24)}$ & $\bar{3}\,2\,\bar{1}\,\bar{4}\,\bar{5}$ & $L_1^{(13)}$ & $\bar{5}\,\bar{2}\,3\,\bar{4}\,\bar{1}$ & $L_1^{(1 5)}$ & $\bar{2}\,\bar{1}\,3\,\bar{5}\,\bar{4}$ & $L_2^{(1 2) (45)}$\\
$1\,\bar{4}\,3\,\bar{2}\,5$ & $L_1^{ (2 4)}$ & $\bar{3}\,2\,\bar{1}\,\bar{5}\,\bar{4}$ & $L_2^{(1 3) (4 5)}$ & $\bar{5}\,\bar{4}\,3\,\bar{2}\,\bar{1}$ & $L_2^{(1 5) (2 4)}$ & $\bar{2}\,\bar{1}\,3\,\bar{4}\,\bar{5}$ & $L_1^{(12)}$\\
$1\,4\,3\,2\,5$ & $L_1^{ (2 4)}$ & $\bar{3}\,2\,\bar{1}\,5\,4$ & $L_2^{(1 3) (4 5)}$ & $\bar{5}\,4\,3\,2\,\bar{1}$ & $L_2^{(1 5) (2 4)}$ & $\bar{2}\,\bar{1}\,\bar{3}\,4\,\bar{5}$ & $L_1^{(12)}$\\
$1\,4\,\bar{3}\,2\,\bar{5}$ & $L_1^{ (2 4)}$ & $3\,2\,1\,5\,4$ & $L_2^{(1 3) (4 5)}$ & $5\,4\,3\,2\,1$ & $L_2^{(1 5) (2 4)}$ & $2\,1\,\bar{3}\,4\,\bar{5}$ & $L_1^{(12)}$\\
$1\,\bar{4}\,\bar{3}\,\bar{2}\,\bar{5}$ & $L_1^{ (2 4)}$ & $3\,2\,1\,\bar{5}\,\bar{4}$ & $L_2^{(1 3) (4 5)}$ & $5\,\bar{4}\,3\,\bar{2}\,1$ & $L_2^{(1 5) (2 4)}$ & $2\,1\,\bar{5}\,4\,\bar{3}$ & $L_2^{(1 2) (3 5)}$\\
$1\,\bar{2}\,\bar{3}\,\bar{4}\,\bar{5}$ & $L_0$ & $3\,2\,1\,\bar{4}\,\bar{5}$ & $L_1^{(13)}$ & $5\,\bar{2}\,3\,\bar{4}\,1$ & $L_1^{(15)}$ & $\bar{2}\,\bar{1}\,\bar{5}\,4\,\bar{3}$ & $L_2^{(1 2) (3 5)}$\\
$1\,\bar{2}\,\bar{4}\,\bar{3}\,\bar{5}$ & $L_1^{(3 4)}$ & $3\,\bar{2}\,1\,4\,\bar{5}$ & $L_1^{(13)}$ & $5\,2\,\bar{3}\,\bar{4}\,1$ & $L_1^{(15)}$ & $\bar{2}\,\bar{1}\,5\,4\,3$ & $L_2^{(1 2) (3 5)}$\\
$1\,\bar{2}\,4\,3\,\bar{5}$ & $L_1^{(3 4)}$ & $\bar{3}\,\bar{2}\,\bar{1}\,4\,\bar{5}$ & $L_1^{(13)}$ & $5\,2\,\bar{4}\,\bar{3}\,1$ & $L_2^{(15) (34)}$ & $2\,1\,5\,4\,3$ & $L_2^{(1 2) (3 5)}$\\
$1\,2\,4\,3\,5$ & $L_1^{(3 4)}$ & $\bar{1}\,\bar{2}\,\bar{3}\,4\,\bar{5}$ & $L_0$ & $5\,2\,4\,3\,1$ & $L_2^{(15) (34)}$ & $2\,1\,3\,4\,5$ & $L_1^{(12)}$\\
\hline
\end{tabular}
\caption{Hamilton cycle with distance 2 for $G(I_5^D)$.}
  \label{table:GI5}
\end{table}

\bigskip

Computational evidence suggests  that the process described above can be generalized for any integer $n\geq 6$, since the availability of connection words that are used  to link two distinct paths $L_i^{\alpha}$ and $L_{i+1}^{(a\,b)\alpha}$ together, satisfying the distance requirement, will increase in number as $n$ gets bigger. This leads to the following conjecture.

\begin{conjecture}
There is an Hamilton cycle in the restriction to involutions of the Cayley graph $G(\ss_n^B,X_1^D\cup X_2^D)$, with Hamming distance two, for $n\geq 4$, where
\begin{align*}
X_1^D&=\{t_{i,j}: i\leq j<n\}\\
X_2^D&=\{t_{i,n}\cdot t_{j,n}:i<j<n\}
\end{align*}
This is the minimal Hamming distance of any Gray code for $I_n^D$.
\end{conjecture}

\end{document}